\documentclass[a4paper,12pt]{article}

\usepackage[pdftex]{graphicx}
\usepackage[section] {placeins}
\usepackage{enumerate}
\usepackage[normalem]{ulem}
\usepackage{amsmath}
\usepackage{amsthm}
\usepackage{latexsym}
\usepackage{amsfonts}
\usepackage{amssymb}

\usepackage{tikz-cd}
\usetikzlibrary{graphs}
\usetikzlibrary{arrows.meta}
\usetikzlibrary{shapes.geometric}

\usepackage{authblk}

\usepackage{latexsym,amssymb,amsfonts, amsthm, amsmath}
\usepackage[english]{babel}
\usepackage{bm}
\usepackage{mathrsfs}

\newtheorem{thm}{Theorem}[section]

\newtheorem{lem}[thm]{Lemma}
\newtheorem{cor}[thm]{Corollary}

\newtheorem{conj}[thm]{Conjecture}

\theoremstyle{definition} 
 
\newtheorem{exa}[thm]{Example}

\newcommand{\Z}{\mathbb{Z}}

\newcommand{\supp}{\mathop{\rm supp}\nolimits}

\newcommand*\boldell{\pmb{\ell}}

\newcommand{\apph}[3]{\mathrel{\operatorname*{{\uplus}^{#1}}_{#2}^{#3}}}

\usepackage{fullpage}

\begin{document}

\title{A Coprime Buratti-Horak-Rosa Conjecture and Grid-Based Linear Realizations}

\author[1]{Onur A\u{g}{\i}rseven}  
\author[2]{M.~A.~Ollis\footnote{Corresponding author.  Email: \texttt{matt\_ollis@emerson.edu}} } 

\affil[1]{Unaffiliated}
\affil[2]{Marlboro Institute for Liberal Arts and Interdisciplinary Studies, Emerson College, Boston, Massachusetts 02116, USA}


\maketitle

\begin{abstract}
We propose a ``Coprime Buratti-Horak-Rosa (BHR) Conjecture": If~$L$ is a multiset of size~$v-1$ with support $\supp(L) \subseteq \{1, 2, \ldots,  \lfloor v/2 \rfloor\}$ such that $\gcd(v,x) = 1$ for all~$x \in L$, then~$L$ is realizable.  This is a specialization of the well-known BHR~Conjecture and it includes Buratti's original conjecture.

We argue that the most effective route to a resolution of the conjecture when $|\supp(L)| = 3$ is to focus on $L = \{1^a, x^b, y^c\}$, where $1<x<y$, with~$a$ large subject to~$a < x+y$.  We use grid-based graphs to construct linear realizations for many such multisets.  A partial list of parameter sets that the constructions cover: 
\begin{itemize}
\item $a = x+y-1$,
\item $a = x+y-2$ when $x=3$ or $x$ is even,
\item  $a \geq 4x-3$ for $x$ odd, $y > 2x-2$, and $b \geq y-2x+2$,
\item $a \geq x$ for $y=tx$, with $x$ and $t$ odd,  and~$b \geq tx+2t-3$, 
\item $a \geq 7$ for~$x=3$ and~$b \geq y-4$.
\end{itemize}  
As well as these (and further) immediate results, the techniques introduced show promise for further development, both to head towards a proof of the conjecture when~$|\supp(L)| = 3$ and for situations with~$| \supp(L) | > 3$.

We also show that if $y > (2x^2 + 2x + 1)/(x-2)$ then the Coprime BHR Conjecture holds for~$\{1^a,x^b,y^c\}$ for infinitely many values of~$v$, and that there are at most~3 values of~$v$ for which it does not hold when~$(x,y) = (6,18)$.

\vspace{3mm}
\end{abstract}

\section{Introduction}\label{sec:intro}

Let~$K_v$ be the complete graph on~$v$ vertices with vertex labels~$\{0, 1, \ldots, v-1\}$.  Define an induced edge labeling, called the {\em length}, by
$$\ell(x,y) = \min ( |y-x|, v - |y-x| ).$$
That is, the length of an edge between vertices~$x$ and~$y$ is the distance between them around a cycle with the~$v$ vertices $\{0,1,\ldots,v-1\}$ in consecutive order. An edge with length~$\ell$ is called an {\em $\ell$-edge}.

Given a subgraph~$G$ of~$K_v$ we obtain a multiset of lengths corresponding to the edges in~$G$.  Say that this multiset is {\em realized} by~$G$ or that~$G$ is a {\em realization} of the multiset.  We are interested in which multisets are realizable by Hamiltonian paths.

This problem was introduced by Buratti~\cite{West07}, who made the conjecture below. (Note: For a multiset~$L$, the set~$\supp(L) = \{ x : x \in L \}$ is the {\em support} of~$L$.)

\begin{conj}\label{conj:buratti}{\rm (Buratti's Conjecture)}
Let~$v$ be prime.  If~$L$ is a multiset of size~$v-1$ with $\supp(L) \subseteq \{1, 2, \ldots,  \lfloor v/2 \rfloor\}$, then~$L$ is realizable.
\end{conj}

This was extended by Horak and Rosa to arbitrary~$v$~\cite{HR09}, with this formulation due to Pasotti and Pellegrini~\cite{PP14b}, and is now commonly known as the Buratti-Horak-Rosa (BHR) Conjecture:

\begin{conj}\label{conj:bhr}
{\rm (BHR Conjecture)}
Let~$L$ be a multiset of size~$v-1$ with $\supp(L) \subseteq \{1, 2, \ldots,  \lfloor v/2 \rfloor\}$. Then~$L$ is realizable if and only if for any divisor~$d$ of~$v$ the number of multiples of~$d$ in~$L$ is at most~$v-d$.
\end{conj}

When~$v$ is prime, the BHR Conjecture reduces to Buratti's Conjecture.  The divisor condition is necessary~\cite{HR09} and the BHR~Conjecture is easily seen to be true when the support has size~1.  

Two early papers independently investigated the case when the support has size~2.  Dinitz and Janisweski showed that Buratti's Conjecture holds in this case and Horak and Rosa, in the same paper that extended Buratti's Conjecture, showed that the BHR~Conjecture holds in this case~\cite{DJ09,HR09}.  Both Dinitz and Janisweski's method and Horak and Rosa's method for prime~$v$ are in fact successful whenever~$v$ is coprime to all elements of the support and this was explicitly noted by Horak and Rosa.

Here we propose an intermediate conjecture that corresponds to this situation.  It reduces to Buratti's Conjecture when~$v$ is prime and includes exactly the instances of the BHR~Conjecture where the divisor condition is vacuous.

\begin{conj}\label{conj:cbhr}{\rm (Coprime BHR Conjecture)}
If~$L$ is a multiset of size~$v-1$ such that $\supp(L) \subseteq \{1, 2, \ldots,  \lfloor v/2 \rfloor\}$ and $\gcd(v,x) = 1$ for all~$x \in L$, then~$L$ is realizable.
\end{conj}

In all three cases, to say that the conjecture holds for a set means that it holds for all multisets with that set as the support.

Much work has been done on multisets with supports larger than~2, but with little resolution (we describe some of this in the next section).  In this paper we focus our attention on the case where the support has size exactly 3.  

For the general situation of a multiset~$L$ with~$L =  \{ x_1^{a_1}, x_2^{a_2}, \ldots, x_k^{a_k}\}$, we may consider the elements to be members of the cyclic group~$\Z_v$ and obtain an equivalent instance of the BHR~Conjecture as follows~\cite{AO1}.  For any $g \in \{1,2, \ldots, v-1\}$, define the {\em reduced form} of~$g$ to be $\widehat{g} = \min(g, v-g)$.  Take~$s$ with~$\gcd(v,s) = 1$, so multiplication by~$s \pmod{v}$ is an automorphism of~$\Z_v$.  A realization of~$L$ is equivalent to a realization of~$\{\widehat{sx_1}^{a_1}, \widehat{sx_2}^{a_2}, \ldots, \widehat{sx_k}^{a_k} \}$.

This has two foundational implications for work on the Coprime BHR Conjecture for supports of size~3.  First, if the multiset in question does not contain~1 then we may translate it into an equivalent multiset that does.  Second, given a multiset of the form~$L = \{1^a, x^b, y^c\}$ with $\gcd(v,x) = \gcd(v,y) = 1$ and~$a+b+c = v-1$, we can obtain two equivalent multisets $L' = \{1^b, \widehat{x^{-1}y}^{\ c}, \widehat{x^{-1}}^{\ a} \}$ and $L'' = \{1^c, \widehat{y^{-1}}^{\ a}, \widehat{xy^{-1}}^{\ b} \}$.  The following notation is sometimes useful:
$$\dot{x} = \min(\  \widehat{x^{-1}y}, \  \widehat{x^{-1}} ) 
\ \ \text{ and } \  \ \dot {y} = \max(\  \widehat{x^{-1}y}, \  \widehat{x^{-1}} ),
$$
\vspace{-6mm}
$$\ddot{x} = \min(\  \widehat{y^{-1}}^{\ a}, \widehat{xy^{-1}}^{\ b} \} ) 
\ \  \text{ and } \ \  \ddot {y} = \max(\  \widehat{y^{-1}}^{\ a}, \widehat{xy^{-1}}^{\ b} \} ),
$$
meaning that $\supp(L') = \{1, \dot{x}, \dot{y} \}$ and $\supp(L'')  = \{1, \ddot{x}, \ddot{y}\}$ with $\dot{x} < \dot{y}$ and $\ddot{x} < \ddot{y}$.

For brevity, call a multiset~$L$ {\em admissible} if it meets the conditions of the BHR Conjecture and {\em strongly admissible} if it also meets the conditions of the Coprime BHR Conjecture.

The starting point of our paper is the following result and the implications that follow from considering equivalent multisets.

\begin{thm}\label{th:1xy}{\rm \cite{AO1}} {\rm (An earlier result)}
Let $L = \{1^a, x^b, y^c\}$ be admissible.  If $a \geq x+y$ then~$L$ is realizable.
\end{thm}

As we shall see in the next section, the result from~\cite{AO1} is slightly stronger, but this is sufficient to illustrate the framework for the paper.

\begin{cor}\label{cor:1xy}{\rm (Counterexample characterization)}
Let $L = \{1^a, x^b, y^c\}$ be strongly admissible with equivalent multisets $L' = \{1^b, \dot{x}^c, \dot{y}^a  \}$ and $L'' = \{1^c, \ddot{x}^a, \ddot{y}^b\}$.   If~$a\geq x+y$, $b \geq \dot{x} + \dot{y}$ or $c \geq \ddot{x} + \ddot{y}$ then~$L$ is realizable.  In particular, any counterexample to the Coprime BHR Conjecture has
$$ x + y + \dot{x} + \dot{y} + \ddot{x} + \ddot{y} \geq v+2. $$
\end{cor}

\begin{proof}
The first statement follows immediately from equivalence and Theorem~\ref{th:1xy}.  The second follows from the first and the fact that~$a + b + c = v-1$.
\end{proof}

Sometimes Corollary~\ref{cor:1xy} is sufficient to resolve all instances of the Coprime BHR Conjecture for a given size and support; usually it is not.

\begin{exa}\label{ex:97a}{\rm (Uses and limitations of the characterization)}
Let~$v=97$.  Suppose~$\supp(L) = \{1,7,13\}$.  Then the equivalent supports are~$\{1,8,15\}$ and~$\{1,12,14\}$.  As
$$7+13+8+15+12+14 = 69 < 99$$
the Coprime BHR Conjecture holds for support~$\{1,7,13\}$ at~$v=97$ (and also for supports~$\{1,8,15\}$ and~$\{1,12,14\}$).

Now suppose~$\supp(L) = \{1,10,13\}$.  The equivalent supports are $\{1,11,29\}$ and $\{1,15,44\}$.  As
$$10+13+11+29+15+44 = 122 \geq 99$$
Corollary~\ref{cor:1xy} does not imply the Coprime BHR Conjecture holds in this case.  
\end{exa}

Although the first case of Example~\ref{ex:97a} is not typical, a heuristic argument suggests that it is also not a vanishingly rare occurrence.  If we make the (false, but empirically reasonable) assumption that~$x,y,\dot{x}, \dot{y},\ddot{x},\ddot{y}$ are six independent random variables taken from the set~$\{2, 3, \ldots, \lfloor v/2 \rfloor \}$ then using a normal approximation to the Irwin-Hall distribution, we obtain a $z$-score of~$- \sqrt{2}$ for the probability that 
$$ x + y + \dot{x} + \dot{y} + \ddot{x} + \ddot{y} < v+2. $$
This corresponds to~7.86\%.  While this is not a rigorous argument, it suggests that existing results already resolve the Coprime BHR Conjecture at any fixed~$v$ for a non-trivial fraction of supports.
However, even the best case scenario, where~7.86\% is trustworthy and on the low side, leaves the remaining 90+\% supports to work on.  

While Corollary~\ref{cor:1xy} is not strong enough to resolve the Coprime BHR Conjecture in the second case of Example~\ref{ex:97a}, it still severely limits what form a counterexample may take.   Notice that we can resolve the conjecture for~$a \geq 23$ and~$b \geq 40$, so any counterexample must have $a+b < 63$; equivalently, $c \geq (97-1) - 63 =  33$.  Thus, looking at the equivalent formulation with support~$\{1,15,44\}$ we can assume that there are at least~33 $1$-edges available.

In the existing work on the conjectures, the most difficult cases seem to be those with a smaller number of 1-edges.  However, the above suggests that the majority of unresolved instances of the Coprime BHR Conjecture (and hence Buratti's Conjecture) are those with a larger number of $1$-edges subject to not being realizable by Theorem~\ref{th:1xy} or its known strengthenings, which are described in the next section.  Moreover, there is a route to a resolution of the conjecture that entirely avoids constructions with very few 1-edges, using the above-mentioned ideas of Dinitz and Janiszewski and of Horak and Rosa~\cite{DJ09,HR09}.

\begin{lem}\label{lem:low1}{\rm (Outstanding cases)}
If the Coprime BHR Conjecture holds for multisets of the form $\{1^a, x^b, y^c\}$ when~$a \geq \max(b,c)$  then it holds for all multisets with support of size~$3$.
\end{lem}

\begin{proof}
Suppose~$L = \{1^a, x^b, y^c\}$.  If~$a \geq \max(b,c)$ then~$L$ is realizable by hypothesis.  Otherwise, if~$b \geq \max(a,c)$ then we may consider the equivalent multiset~$L'$, which has one of the forms $\{1^b, \dot{x}^c, \dot{y}^a  \}$ or~$\{1^b, \dot{x}^a, \dot{y}^c \}$ and is realizable by hypothesis.  A similar argument applies with~$L''$ if~$c \geq \max(a,b)$.
 \end{proof}

An analogous argument may be applied for multisets with larger support.  

It is unlikely that Lemma~\ref{lem:low1} is part of the most straightforward resolution of the Coprime BHR Conjecture: sometimes realizations can be constructed for smaller values of~$a$ with relative ease and these can remove the need to directly solve difficult cases with~$a \geq \max(b,c)$. We revisit this at the end of the Section~\ref{sec:prev} when we have more tools and background at our disposal; see Theorem~\ref{th:all_lin} and the subsequent discussion.

\bigskip

These observations drive our choices of constructions.   When working on this problem, it is all too easy to be drawn into ever more elaborate constructions that give ever decreasing coverage, especially as this is sometimes necessary to close out a particular family of cases.  However, our goal in this paper is to present constructions chosen for their (subjectively judged) simplicity to effectiveness ratio in exploring the ``target region" of parameters rich in unresolved instances, i.e. those for strongly admissible~$\{1^a, x^b, y^c\}$ and~$a$ large subject to~$a < x+y$. Our results reflect that goal.

\begin{thm}\label{th:main}{\rm (Select Main Results)}
Let $L = \{1^a, x^b, y^c \}$ with $1< x<y$.  Then~$L$ has a linear realization in the following cases:
\begin{enumerate}
\item $a \geq x+y-1$,
\item $a \geq x+y-2$ when $x=3$ or $x$ is even,
\item $a \geq 4x-3$ for $x$ odd, $y > 2x-2$, and $b \geq y-2x+2$,
\item $a \geq x$ for $y=tx$, with $x$ and $t$ odd,  and~$b \geq tx+2t-3$, 
\item $a \geq 7$ for~$x=3$ and~$b \geq y-4$.
\end{enumerate}
\end{thm}

Section~\ref{sec:prev} introduces various useful ideas from previous considerations of the BHR Conjecture and includes a full account of the current state of knowledge regarding realizations with supports of size~3.  Sections~\ref{sec:red} and~\ref{sec:1xtx} develop the bulk of the new constructions.  The main idea is to start with an appropriate optimal construction for a multiset with support of size~2 and replace some of the edges with those of the desired third type.  

In Section~\ref{sec:red}, we switch from supports of the form~$\{1,y\}$ to~$\{1,x,y\}$ by replacing 1-edges with $x$-edges, which is especially effective when~$y$ is considerably larger than~$x$.  Theorem~\ref{th:1xy_new} improves Corollary~\ref{cor:1xy}, proving the first main result, while Theorem~\ref{th:odd_hops} proves the third.  Theorem~\ref{th:xoddyeven}, on the other hand, is useful when $x$ is closer to $y$ in size.

In Section~\ref{sec:1xtx}, we switch from supports of the form~$\{1,x\}$  to~$\{1,x,tx\}$ by replacing 1-edges with $tx$-edges.  Theorem~\ref{th:1xtx_gen} provides a bound on $a$ for odd $x$ that is very strong compared to most known results.  Theorem~\ref{th:1xtx_tight} improves this further. In all but one of its cases, including the fourth main result, the bound on~$b$ is a constant depending only on~$t$ and~$x$, while in the remaining case, the bound on~$b$ as $v$ increases is improved by a factor of~$t$, compared to Theorem~\ref{th:1xtx_gen}.

Section~\ref{sec:conseq} gives significant improvements to existing results for the case~$L = \{1^a, 3^b, y^c\}$ and provides the remaining main results.  It includes a detailed examination of how our constructions apply to the support~$\{1,6,18\}$, reducing the possibility of a counterexample to the Coprime BHR Conjecture with this support to having~$v = 47, 49$ or~$59$.

Sections~\ref{sec:1xtx} and~\ref{sec:conseq} also consider various consequences for the conjectures. The ones included below follow from Theorems~\ref{th:1xtx},~\ref{th:13ynum}, and~\ref{th:bigy}.

\begin{thm}\label{th:bhr}{\rm (Select Implications for BHR)}
For a multiset~$L$ of size~$v-1$ with $\supp(L) = \{1,x,y\}$, the BHR Conjecture holds in the following cases:
\begin{enumerate}
\item $x>2$, $y > 2x + 19$, $v \equiv \pm 1 \pmod{xy}$, and $v$ is sufficiently large,
\item $x=3$, $y=3t$, $v \equiv \pm 1 \pmod{3t}$, and $v > 6t+53$,
\item $x,t \geq 7$, $y=tx$, ~$v \equiv \pm 1 \pmod{tx}$.
\end{enumerate}
\end{thm}

\section{Notation, Tools, and Previous Work}\label{sec:prev}

As mentioned in the introduction, it is known that we can do better than Theorem~\ref{th:1xy}'s result that realizations for~$L = \{1^a, x^b, y^c\}$ exist whenever~$a \geq x+y$.  Theorem~\ref{th:1xy2} gives general results from~\cite{AO1,OPPS} and Theorem~\ref{th:known} gives special cases from the literature that do better still under certain conditions.

\begin{thm}\label{th:1xy2}{\rm \cite{AO1,OPPS}} {\rm (Known results for arbitrary supports of size 3)}
Let~$L = \{1^a, x^b, y^c\}$, with~$1<x<y$, be a multiset of size~$v-1$.  Then~$L$ is realizable when:
\begin{enumerate}
    \item $x,y$ odd and~$a \geq x+y$,
    \item $x$ odd, $y$ even and~$a \geq x+y-1$,
    \item $x$ even, $y$ odd and $a \geq x+y-2$,
    \item $x,y$ even and $a \geq y-1$
\end{enumerate}
\end{thm}

\begin{thm}\label{th:known}{\rm (Known results for specific types of supports of size 3)}
Let $L$ be a multiset of size~$v-1$ with~$|\supp(L)| = 3$. If $L$ is admissible then it is realizable when:
\begin{enumerate}
\item $v \leq 37$ \cite{MP,Meszka},
\item 
  $\max(L) \leq 7$ or $\supp(L) = \{ 1,2,8 \}, \{1,2,10\}, \{1,2,12 \}$  \cite{CD10,CO,PP14, PP14b}, 
\item $L = \{1^a,2^b, y^c  \}$ when $a+b \geq y-1$ and either~$y$ is even or $a \geq 3$ \cite{AO1,PP14},
\item $L = \{1^a, 3^b, y^c \}$ when~$y$ is even, $c$ is odd, and either~$a \geq y+1$ or $a=y$ and $3 \nmid b$~\cite{Avila23}.
\item $L = \{1^a, x^b, (x+1)^c\}$ and~$a \geq x+1$ \cite{AO1},
\item $L = \{ 1^a, x^b, (2x)^c \}$ when $a \geq x-2$, $c$ is even and $b \geq 5x-2+c/2$ \cite{OPPS2}.
\item \label{th:evenx} $L = \{1^a, x^b, y^c\}$ with $x$ even and $y\geq 2x+1$, when $a \geq 3x-2$ and  $a+b \geq x+y-1$ \cite{AO1}.
\end{enumerate}
\end{thm}

Further, in~\cite{AO1}, arguments based on the number theoretic approach outlined in the introduction are used to deduce the following result from the third and fifth items of Theorem~\ref{th:known}.

\begin{thm}\label{th:ao_largev}{\rm \cite{AO1}} {\rm (BHR implications deduced)}
The Coprime BHR Conjecture holds in the following cases:
\begin{enumerate}
    \item $L = \{1^a,2^b,y^c\}$ and $v > 4y$, except possibly when~$y$ is odd and~$a \in \{1,2\}$,
    \item $\supp(L) = \{ 1,x,x+1\}$ and $v \geq 2x^2 + 13x + 11$.
\end{enumerate}
\end{thm}

Most of the results in Theorem~\ref{th:known} were proved via constructions of special types of realizations.  We introduce some of these types now.

The {\em linear length} of an edge between vertices~$x$ and~$y$ in~$K_v$ is defined to be~$|x-y|$.  We can take the multiset~$L$ of linear lengths of edges of a subgraph~$G$; in this case say that~$L$ is {\em linearly realized} by~$G$ or that~$G$ is a {\em linear realization} of~$L$.  We are again interested in~$G$ being a Hamiltonian path.  (We sometimes use ``cyclic realization" for realization as previously defined when it clarifies the distinction.)

If a linear realization has linear lengths at most~$\lfloor v/2 \rfloor$ then it is a cyclic realization for the same multiset.  More generally, any linear realization for a multiset~$L$ is a cyclic realization for~$\{\widehat{x} : x \in L\}$.  Given this close relationship it usually does not cause confusion to refer to linear lengths simply as lengths; we shall do so and disambiguate where necessary.
 
A realization is {\em standard} if its first vertex is~0.  It is {\em perfect} if it is standard and also the last vertex is~$v-1$.

The {\em complement} of a realization is formed by replaced each vertex~$x$ with~$v-1-x$.  It realizes the same multiset.  The {\em (embedded) translation} of a realization by~$t$ is formed by adding~$t$ to each vertex.  It realizes the same multiset, as a non-Hamiltonian path in~$K_{v+t}$.

\begin{lem}\label{lem:concat}{\rm \cite{HR09,OPPS}} {\rm (Concatenation)}
If~$L$ and~$M$ have standard linear realizations, then~$L \cup M$ has a linear realization.  If the realization for~$M$ is perfect, then~$L \cup M$ has a standard linear realization.  If both are perfect, then~$L \cup M$ has a perfect linear realization.

If~$L$ has a standard (respectively perfect) linear realization then for all~$s \geq 0$ the multiset~$L \cup \{1^s\}$ has a standard (respectively perfect) linear realization.
\end{lem}

\begin{proof}[Proof Construction.]
Suppose~$L$ has size~$v-1$ and~$\mathbf{g}$ is the standard linear realization for~$L$.   Suppose~$M$ has size~$w-1$ and~$\mathbf{h}$ is the standard linear realization for~$L$.   Define the {\em concatenation}~$\mathbf{g} \oplus \mathbf{h}$ to be the Hamiltonian path in~$K_{v+w-1}$ obtained by taking complement of~$\mathbf{g}$ and the translation of~$\mathbf{h}$ by~$v-1$ and identifying the vertices labeled~$v-1$.  This concatenation has the required properties for the first paragraph of the statement.

For the second statement, use the same concatenation with~$\mathbf{h} = [ 0, 1, \ldots, s ] $ as the perfect linear realization for~$\{1^s\}$. 
\end{proof}

In~\cite{AO1}, grid-based graphs, a new tool for visualizing and hence thinking about certain types of realizations, were introduced.  The same method is at the heart of our constructions here.  

\begin{figure}[tp]
\caption{Standard linear realizations for~$\{1^7, 8^{22}  \}$ and~$\{1^5, 5^{2}, 8^{21}   \}$.}\label{fig:grid_eg}
\begin{center}
\begin{tikzpicture}[scale=0.9, every node/.style={transform shape}]
\fill (0,1) circle (2pt) ;\fill (0,2) circle (2pt) ;\fill (0,3) circle (2pt) ;\fill (0,4) circle (2pt) ;\fill (1,1) circle (2pt) ;\fill (1,2) circle (2pt) ;\fill (1,3) circle (2pt) ;\fill (1,4) circle (2pt) ;\fill (2,1) circle (2pt) ;\fill (2,2) circle (2pt) ;\fill (2,3) circle (2pt) ;\fill (2,4) circle (2pt) ;\fill (3,1) circle (2pt) ;\fill (3,2) circle (2pt) ;\fill (3,3) circle (2pt) ;\fill (3,4) circle (2pt) ;\fill (4,1) circle (2pt) ;\fill (4,2) circle (2pt) ;\fill (4,3) circle (2pt) ;\fill (4,4) circle (2pt) ;\fill (5,1) circle (2pt) ;\fill (5,2) circle (2pt) ;\fill (5,3) circle (2pt) ;\fill (5,4) circle (2pt) ;\fill (6,1) circle (2pt) ;\fill (6,2) circle (2pt) ;\fill (6,3) circle (2pt) ;\fill (7,1) circle (2pt) ;\fill (7,2) circle (2pt) ;\fill (7,3) circle (2pt) ;
\draw (0,1) -- (0,4) -- (1,4) -- (1,1) -- (2,1) -- (2,4) -- (3,4) -- (3,1) -- (4,1) -- (4,4) 
   -- (5,4) -- (5,1) -- (6,1) -- (6,3) -- (7,3) -- (7,1)  ;
\node at (-0.3, 1) {\tiny 0} ;  \node at (-0.3, 2) {\tiny 8} ;  \node at (-0.3, 4) {\tiny 24} ; \node at (5.3, 4) {\tiny 29} ; \node at (7.3, 1) {\tiny 7} ;  \node at (7.3, 3) {\tiny 23} ;   \fill (9,1) circle (2pt) ;\fill (9,2) circle (2pt) ;\fill (9,3) circle (2pt) ;\fill (9,4) circle (2pt) ;\fill (10,1) circle (2pt) ;\fill (10,2) circle (2pt) ;\fill (10,3) circle (2pt) ;\fill (10,4) circle (2pt) ;\fill (11,1) circle (2pt) ;\fill (11,2) circle (2pt) ;\fill (11,3) circle (2pt) ;\fill (11,4) circle (2pt) ;\fill (12,1) circle (2pt) ;\fill (12,2) circle (2pt) ;\fill (12,3) circle (2pt) ;\fill (12,4) circle (2pt) ;\fill (13,1) circle (2pt) ;\fill (13,2) circle (2pt) ;\fill (13,3) circle (2pt) ;\fill (14,1) circle (2pt) ;\fill (14,2) circle (2pt) ;\fill (14,3) circle (2pt) ;\fill (15,1) circle (2pt) ;\fill (15,2) circle (2pt) ;\fill (15,3) circle (2pt) ; \fill (16,0) circle (2pt) ;\fill (16,1) circle (2pt) ; \fill (16,2) circle (2pt) ; \fill (16,3) circle (2pt) ;
\draw (9,1) -- (9,4) -- (10,4) -- (10,1); \draw (11,1) -- (11,4) -- (12,4) -- (12,1) -- (13,1) -- (13,3) -- (14,3) -- (14,1); \draw (15,1) -- (15,3)--(16,3)--(16,0) ; \draw  plot [smooth] coordinates {(9,1) (11.5,0.3) (14,1)}; \draw  plot [smooth] coordinates {(10,1) (12.5,0.3) (15,1)};
\node at (8.7, 1) {\tiny 1} ;  \node at (8.7, 2) {\tiny 9} ;  \node at (8.7, 4) {\tiny 25} ;  \node at (11.3, 1) {\tiny 3} ;\node at (12.3, 4) {\tiny 28} ;  \node at (16.3, 0) {\tiny 0} ;  \node at (16.3, 1) {\tiny 8} ;   \node at (16.3, 3) {\tiny 24} ;  
\end{tikzpicture}\end{center}
\end{figure}
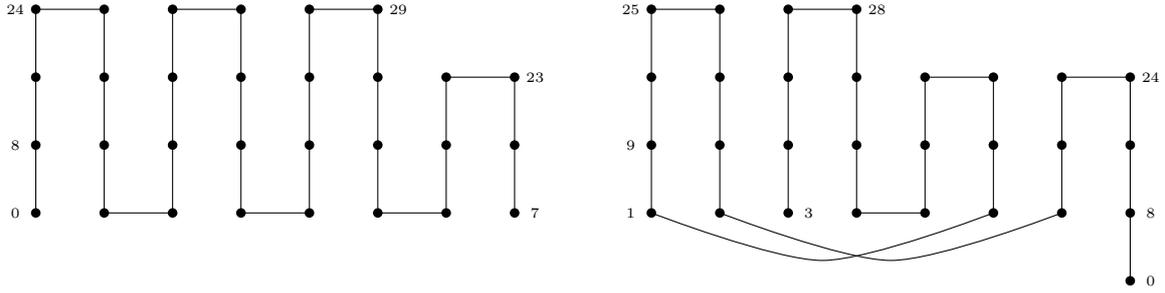

 We place the vertices of~$K_v$ at $v$ points of the integer lattice, labeled in increasing order going from left to right and then bottom to top.  We arrange them so that there are (usually)~$x$ vertices in a row, for some~$x \in L$, meaning that horizontal edges between adjacent vertices are $1$-edges and vertical ones are $x$-edges.  Other edge lengths in this paper also have fixed geometry; there will only be one other edge type at a time.  Figure~\ref{fig:grid_eg} illustrates diagrams for realizations of~$\{1^7, 8^{22}  \}$ and~$\{1^5, 5^{2}, 8^{21}   \}$. 

Notice that the realizations in Figure~\ref{fig:grid_eg} use all of the elements that are congruent to some~$k \pmod{x}$ (where an $x$-edge is between vertically adjacent vertices and in both examples~$x=8$) consecutively.  This is common in the constructions in this paper and all of them tend to use most such elements consecutively.  Following~\cite{AO1}, we use this to motivate some definitions and notation.

Let~$v= qx+r$ with~$0 \leq r < q$.  Define the {\em fauxset} $\phi_k$ to be the subset of elements in~$\supp(L)$ that are congruent to~$k \pmod{x}$ and let the {\em fauxset traversal}~$\Psi_k$ be the path that traverses the fauxset $\phi_k$ in sequence, realizing only $x$-edges.

More generally, set
$$q^* = q^*(v,r,k) =  
\begin{cases} 
q & \mathrm{if \  } k<r, \\
 q-1 & \mathrm{otherwise.}
\end{cases}
$$
That is, $q^*$ is the value that makes $q^*x + k$ the largest element of~$\varphi_k$.
 Let the {\em fauxset segment} $\Psi_k^{(a,b)}$ be the path 
$[ ax + k, (a+1)x+k, \ldots, bx+k ]$ of elements in~$\varphi_k$ that realizes only $x$-edges.  It is sometimes useful to refer to a path that is the single vertex~$ax+k \in \phi_k$ using this notation; we use~$\Psi_k^{(a)}$.

We also need to connect fauxsets and other paths.  Given paths~$\mathbf{p_1}$ with~$g$ as an end-point and~$\mathbf{p_2}$ with~$h$ has an end-point.  If the length (whether linear or cyclic) of the edge between~$g$ and~$h$ is~$\ell$ then we call the path obtained by adding the edge between~$g$ and~$h$ a {\em bridge concatenation} and denote it as~$\mathbf{p_1} \uplus^\ell \mathbf{p_2}$.  Although this does not necessarily define a unique construction, it will always be clear from context which option to use when there is a choice.

More generally, given paths~$\mathbf{p_1}, \mathbf{p_2}, \ldots, \mathbf{p_t}$ and a vector of lengths~$\boldell  = (\ell_1, \ell_2, \ldots, \ell_{t-1})$,  let
$$ \apph{\boldell}{k=1}{t}  \mathbf{p_k}  = 
\mathbf{p_1} \uplus^{\ell_1} \mathbf{p_2} \uplus^{\ell_2} \cdots \uplus^{\ell_{t-1}} \mathbf{p_t}
{\rm \ \ \ \  and \ \ \ \  } 
\apph{\boldell}{k=t}{1}  \mathbf{p_k} = 
\mathbf{p_t} \uplus^{\ell_1} \mathbf{p_{t-1}} \uplus^{\ell_2} \cdots \uplus^{\ell_{t-1}} \mathbf{p_1}.$$
When~$\boldell = (\ell,\ell,\ldots, \ell)$ we use~$\ell$ in place of~$\boldell$ and if~$\ell=1$ we omit it.  With this notation, the realizations in Figure~\ref{fig:grid_eg} are
$$ \apph{}{k=0}{7} \Psi_k  {\rm \ \  \ and \  \ \ } 
\apph{\boldell}{k=8}{1}  \Psi_{h_k}$$
with~$\mathbf{h} = [h_1, \ldots, h_8] = [0,7,2,1,6,5,4,3]$ whose differences are~$\boldell = (1,5,5,1,1,1,1)$.  

To move from one fauxset to another (with respect to~$x$) in a linear realization, it is necessary to use an edge that has a length that is not a multiple of~$x$.  This gives a general bound (see \cite[Lemma~2.1]{AO1}) from which we obtain the following result that is most relevant for our situation:

\begin{lem}\label{lem:hop}{\rm (Connecting fauxsets)}
Let~$L = \{1^a, x^b, y^c \}$ with~$1< x< y$ and~$c>0$.  Any linear realization for~$L$ has~$a + b \geq y-1$.
\end{lem}

The bound in Lemma~\ref{lem:hop} cannot always be obtained. 

Consider the case~$b=0$ and the extra condition that we require a standard realization.  We have~$L = \{1^a,y^c\}$.  Let~$\omega(y,c)$ denote the smallest value of~$a$ such that a standard realization of~$L$ exists (and hence by Lemma~\ref{lem:concat} a standard realization also exists for all larger values of~$a$).  

\begin{thm}\label{th:omega}{\rm \cite{AO1}} {\rm ($\omega$-constructions)}
Let~$y\geq2$ and~$c > 0$.  Write~$c = q'y + r'$ with $0 \leq r' < y$.  If $q'=0$, $r'=1$ or at least one of $y$ and~$r'$ is even then~$\omega(y,c) = y-1$.  Otherwise, $\omega(y,c) = y$.
\end{thm}

\begin{proof}[Proof Construction.]
For the cases with~$\omega(y,c) = y-1$, at least one of the following two constructions gives a successful realization:
$$\mathbf{h_1} = \apph{}{k=0}{y-1} \Psi_k  {\rm \ \ \ and \ \ \ }  \mathbf{h_2} =  \Psi_0 \uplus \left( \apph{}{k=y-1}{1} \Psi_k \right) . $$
The sequence~$\mathbf{h_1}$ is successful when~$r'$ is even; $\mathbf{h_2}$ is successful when~$r'=1$, when~$q'=0$, and when~$y$ and~$r'$ have different parities (which includes the remaining uncovered case with~$y$ even and~$r'$ odd).

When none of the conditions for meeting~$\omega(y,c) = y-1$ are met, set~$q = q'+1$ and $r = r'+1$, so $v = qy + r$, and let
$$\mathbf{h_3} = \Psi_0 \uplus \left( \apph{}{k=y-1}{3} \Psi_k \right) \uplus \Psi_2^{(q^*)} \uplus \Psi_1 \uplus \Psi_2^{(0, q^*-1)}. $$ 
\end{proof}

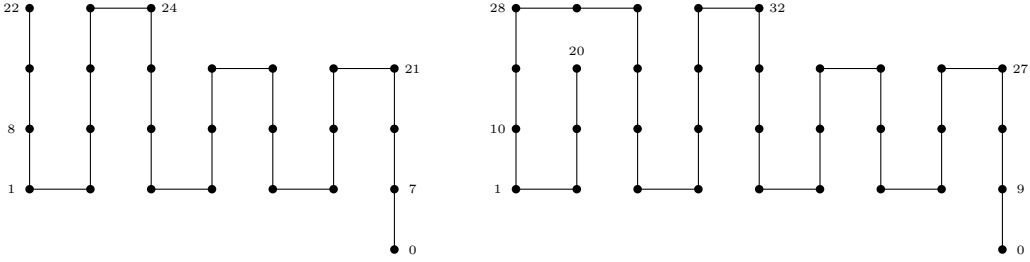
\begin{figure}[tp]
\caption{The standard linear realization $\mathbf{h_2}$ for~$\{1^6, 7^{18}\}$ and the standard linear realization~$\mathbf{h_3}$ for~$\{1^9, 9^{23}\}$.}\label{fig:omega_eg}
\begin{center}
\begin{tikzpicture}[scale=0.8, every node/.style={transform shape}]

\fill (0,1) circle (2pt) ;\fill (0,2) circle (2pt) ;\fill (0,3) circle (2pt) ;\fill (0,4) circle (2pt) ;\fill (1,1) circle (2pt) ;\fill (1,2) circle (2pt) ;\fill (1,3) circle (2pt) ;\fill (1,4) circle (2pt) ;\fill (2,1) circle (2pt) ;\fill (2,2) circle (2pt) ;\fill (2,3) circle (2pt) ;\fill (2,4) circle (2pt) ;\fill (3,1) circle (2pt) ;\fill (3,2) circle (2pt) ;\fill (3,3) circle (2pt) ;\fill (4,1) circle (2pt) ;\fill (4,2) circle (2pt) ;\fill (4,3) circle (2pt) ;\fill (5,1) circle (2pt) ;\fill (5,2) circle (2pt) ;\fill (5,3) circle (2pt) ;\fill (6,0) circle (2pt) ;\fill (6,1) circle (2pt) ;\fill (6,2) circle (2pt) ;\fill (6,3) circle (2pt) ;
\draw (0,4) -- (0,1) -- (1,1) -- (1,4) -- (2,4) -- (2,1) -- (3,1)  -- (3,3) -- (4,3) -- (4,1) -- (5,1) -- (5,3)  -- (6,3) -- (6,0) ;
\node at (-0.3, 1) {\tiny 1} ;  \node at (-0.3, 2) {\tiny 8} ;  \node at (-0.3, 4) {\tiny 22} ; \node at (6.3, 0) {\tiny 0} ;    \node at (6.3, 1) {\tiny 7} ;  \node at (6.3, 3) {\tiny 21} ;   \node at (2.3, 4) {\tiny 24} ;  
\fill (8,1) circle (2pt) ;\fill (8,2) circle (2pt) ;\fill (8,3) circle (2pt) ;\fill (8,4) circle (2pt) ;\fill (9,1) circle (2pt) ;\fill (9,2) circle (2pt) ;\fill (9,3) circle (2pt) ;\fill (9,4) circle (2pt) ;\fill (10,1) circle (2pt) ;\fill (10,2) circle (2pt) ;\fill (10,3) circle (2pt) ;\fill (10,4) circle (2pt) ;\fill (11,1) circle (2pt) ;\fill (11,2) circle (2pt) ;\fill (11,3) circle (2pt) ;\fill (11,4) circle (2pt) ;\fill (12,1) circle (2pt) ;\fill (12,2) circle (2pt) ;\fill (12,3) circle (2pt) ;\fill (12,4) circle (2pt) ;\fill (13,1) circle (2pt) ;\fill (13,2) circle (2pt) ;\fill (13,3) circle (2pt) ;\fill (14,1) circle (2pt) ;\fill (14,2) circle (2pt) ;\fill (14,3) circle (2pt) ;\fill (15,1) circle (2pt) ;\fill (15,2) circle (2pt) ;\fill (15,3) circle (2pt) ;\fill (16,0) circle (2pt) ;\fill (16,1) circle (2pt) ;\fill (16,2) circle (2pt) ;\fill (16,3) circle (2pt) ;
\draw (9,3) -- (9,1) -- (8,1)   -- (8,4) -- (10,4) -- (10,1) -- (11,1) -- (11,4) -- (12,4) -- (12,1) -- (13,1) -- (13,3) -- (14,3) -- (14,1) -- (15,1)  -- (15,3) -- (16,3) -- (16,0) ;
\node at (7.7, 1) {\tiny 1} ;  \node at (7.7, 2) {\tiny 10} ;  \node at (7.7, 4) {\tiny 28} ;  \node at (9, 3.3) {\tiny 20} ;\node at (12.3, 4) {\tiny 32} ;  \node at (16.3, 0) {\tiny 0} ;  \node at (16.3, 1) {\tiny 9} ;  \node at (16.3, 3) {\tiny 27} ;  

\end{tikzpicture}\end{center}
\end{figure}

We use the realizations~$\mathbf{h_1}$,~$\mathbf{h_2}$ and~$\mathbf{h_3}$ from the proof of Theorem~\ref{th:omega} at various points in the paper.  Following~\cite{AO1}, we call them {\em $\omega$-constructions}.
The realization for~$\{1^7,8^{22}\}$ in Figure~\ref{fig:grid_eg} is an example of~$\mathbf{h_1}$. Figure~\ref{fig:omega_eg} gives examples of~$\mathbf{h_2}$ and~$\mathbf{h_3}$. We refer to the type of move transforming~$\mathbf{h_2}$ into~$\mathbf{h_3}$ for given $x$ and $v$, informally, as a {\em tail curl}.

Note that Theorem~\ref{th:omega} and Lemma~\ref{lem:concat} are almost sufficient to prove the first three cases of Theorem~\ref{th:1xy2} (the case with~$x$ even and~$y$ odd needs more work, see~\cite[Lemma~5.1]{AO1}).  The main approach of this paper is to follow this method, but either to replace some of the 1-edges in an $\omega$-construction with support~$\{1, y\}$ with $x$-edges or to replace some of the~$x$-edges in an $\omega$-construction with support~$\{1,x\}$ with~$y$-edges, thus hitting parameters in the target region not covered by Theorem~\ref{th:1xy2}.

\bigskip

To conclude this section we consider the possibility of proving the Coprime~BHR Conjecture using only linear realizations in light of the constraint imposed by Lemma~\ref{lem:hop}.

\begin{thm}\label{th:all_lin}{\rm (An alternative target region)}
If the Coprime BHR Conjecture holds for all strongly admissible multisets of the form $\{1^a, x^b, y^c\}$ with~$x<y$ and $a+b \geq y-1$ then it holds for all multisets with support of size~$3$.
\end{thm}

\begin{proof}
Let~$L = \{1^a, x^b, y^c\}$ with $x<y$ and suppose that $L' = \{1^b, \dot{x}^c, \dot{y}^a  \}$ and $L'' = \{1^c, \ddot{x}^a, \ddot{y}^b\}$ are equivalent formulations of an instance of the Coprime BHR Conjecture. 

By hypothesis, we can realize these multisets if $a+b \geq y-1$, $b+c \geq \dot{y}-1$ or $a+c \geq \ddot{y}-1$.  Therefore any counterexample has
$$ (a+b) + (b+c) + (a+c) < (y-1) + (\dot{y}-1) + (\ddot{y}-1).$$
Using $a+b+c = v-1$, we find the equivalent inequality
$$ 2v + 1 <  y+\dot{y}+\ddot{y}.$$
But each of $y$, $\dot{y}$ and $\ddot{y}$ are at most $v/2$, giving the contradiction $2v+1 < 3v/2$.
\end{proof}

Compared with Lemma~\ref{lem:low1}, we have the disadvantage of needing to cover cases where the number of 1-edges is very small.  However, the advantage of having linear realizations potentially cover most of the parameters is significant, especially given the more advanced state of the tools for this situation.

Although the hypothesis of Theorem~\ref{th:all_lin} as stated cannot be satisfied with linear realizations alone, its result can be generalized by modifying the proof, where the condition ``$a+b \geq y-1$" is replaced with ``$a+b \geq y + \epsilon$" for any fixed $\epsilon \geq 0$ and the result applies for sufficiently large~$v$ (how large is sufficient depends only on $\epsilon$; approximately $v > 6\epsilon$).  This permits the possibility of linear realizations contributing even more to a proof of the Coprime BHR Conjecture for multisets with support of size~3.

We believe that the ultimate route to successful resolution of the Coprime BHR Conjecture for supports of size~3 will use a combination of the approaches in Theorem~\ref{th:all_lin} and Lemma~\ref{lem:low1}.  The remainder of the paper is concerned with extending the inventory of multisets for which we can construct a linear realization (and of multisets that are equivalent to such multisets) as progress towards such a resolution.

\section{Replacing 1-edges with $x$-edges}\label{sec:red}

In this section, the goal is to construct some standard linear realizations for $\{1^a, x^b, y^c\}$ with~$1<x<y$ that have $a < y-1$ and small~$b$.  We mostly limit our attention to odd~$x$.  Similar, but weaker, results are possible for even~$x$, but they are rarely stronger than both Theorem~\ref{th:1xy2}.4 and Theorem~\ref{th:known}.7.

\begin{figure}[tp]
\caption{Standard linear realizations for~$\{1^2, 3^4, 7^{15}\}$ and~$\{1^3, 3^5, 9^{22}\}$.}\label{fig:fullfauxset}
\begin{center} \begin{tikzpicture}[scale=0.8, every node/.style={transform shape}]
\fill (0,1) circle (2pt) ; \fill (0,2) circle (2pt) ; \fill (0,3) circle (2pt) ; \fill (1,1) circle (2pt) ; \fill (1,2) circle (2pt) ; \fill (1,3) circle (2pt) ; \fill (2,1) circle (2pt) ; \fill (2,2) circle (2pt) ;\fill (2,3) circle (2pt) ; \fill (3,1) circle (2pt) ; \fill (3,2) circle (2pt) ; \fill (3,3) circle (2pt) ; \fill (4,1) circle (2pt) ; \fill (4,2) circle (2pt) ; \fill (4,3) circle (2pt) ; \fill (5,1) circle (2pt) ; \fill (5,2) circle (2pt) ; \fill (5,3) circle (2pt) ; \fill (6,0) circle (2pt) ; \fill (6,1) circle (2pt) ; \fill (6,2) circle (2pt) ; \fill (6,3) circle (2pt) ;
\draw (2,3) -- (2,1) ;  \draw (5,1) -- (5,3) -- (4,3) -- (4,1) ; \draw (1,1)  -- (1,3) -- (0,3) -- (0,1) ; \draw (3,1) -- (3,3)  ; \draw (6,3) -- (6,0) ; \draw  plot [smooth] coordinates {(0,1) (1.5,0.5) (3,1)}; \draw  plot [smooth] coordinates {(1,1) (2.5,0.5) (4,1)}; \draw  plot [smooth] coordinates {(2,1) (3.5,0.5) (5,1)}; \draw  plot [smooth] coordinates {(3,3) (4.5,3.5) (6,3)};
\node at (-0.3, 1) {\tiny 1} ;  \node at (-0.3, 2) {\tiny 8} ;  \node at (-0.3, 3) {\tiny 15} ; \node at (6.3, 0) {\tiny 0} ;  \node at (6.3, 1) {\tiny 7} ;  \node at (6.3, 3) {\tiny 21} ;  \node at (2, 3.3) {\tiny 17} ;  
\fill (8,1) circle (2pt) ; \fill (8,2) circle (2pt) ; \fill (8,3) circle (2pt) ; \fill (8,4) circle (2pt) ; \fill (9,1) circle (2pt) ; \fill (9,2) circle (2pt) ; \fill (9,3) circle (2pt) ; \fill (9,4) circle (2pt) ; \fill (10,1) circle (2pt) ; \fill (10,2) circle (2pt) ; \fill (10,3) circle (2pt) ; \fill (10,4) circle (2pt) ; \fill (11,1) circle (2pt) ; \fill (11,2) circle (2pt) ; \fill (11,3) circle (2pt) ; \fill (11,4) circle (2pt) ; \fill (12,1) circle (2pt) ; \fill (12,2) circle (2pt) ; \fill (12,3) circle (2pt) ; \fill (13,1) circle (2pt) ; \fill (13,2) circle (2pt) ; \fill (13,3) circle (2pt) ; \fill (14,1) circle (2pt) ; \fill (14,2) circle (2pt) ; \fill (14,3) circle (2pt) ; \fill (15,1) circle (2pt) ; \fill (15,2) circle (2pt) ;\fill (15,3) circle (2pt) ; \fill (16,1) circle (2pt) ; \fill (16,2) circle (2pt) ; \fill (16,3) circle (2pt) ;
\draw (9,1) -- (9,4) -- (10,4) -- (10,1);\draw (8,1)   -- (8,4)  ;\draw (11,1) -- (11,4) ;\draw (12,3) -- (12,1) ;\draw (13,1) -- (13,3) ;\draw (14,1) -- (14,3) -- (15,3)  -- (15,1) -- (16,1) -- (16,3) ;
\draw  plot [smooth] coordinates {(8,4) (9.5,4.5) (11,4)};\draw  plot [smooth] coordinates {(13,3) (14.5,3.5) (16,3)};\draw  plot [smooth] coordinates {(9,1) (10.5,0.5) (12,1)};\draw  plot [smooth] coordinates {(10,1) (11.5,0.5) (13,1)};\draw  plot [smooth] coordinates {(11,1) (12.5,0.5) (14,1)};
\node at (7.7, 1) {\tiny 0} ;  \node at (7.7, 2) {\tiny 9} ;  \node at (7.7, 4) {\tiny 27} ;  \node at (12, 3.3) {\tiny 22} ;\node at (11.3, 4) {\tiny 30} ; \node at (16.3, 1) {\tiny 8} ;  \node at (16.3, 3) {\tiny 26} ;  
\end{tikzpicture}\end{center}
\end{figure}
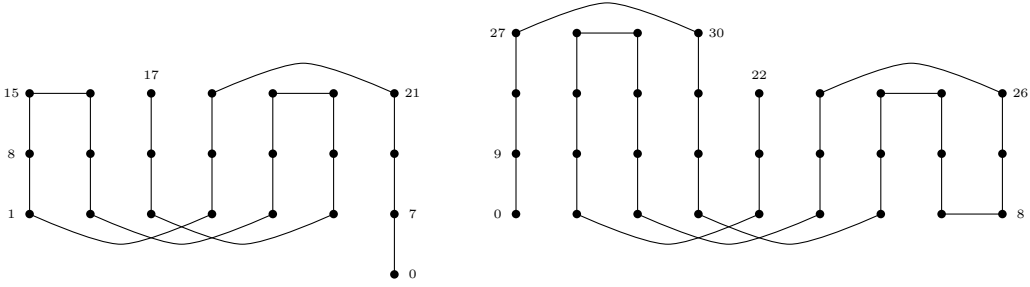

The method is to take the appropriate $\omega$-construction and replacing some of the 1-edges with $x$-edges.  Figure~\ref{fig:fullfauxset} illustrates the type of realization that results.  
These standard linear realizations may then be optionally concatenated with a perfect realization with support~$\{ 1 \}$ and then an $\omega$-construction with support~$\{1,x\}$ to produce realizations with support~$\{1,x,y\}$ that land in the target region.

We start with a result for the cleanest situations: when all fauxsets have equal size or there is a single one larger than the others.  This can be viewed as a special case of later methods although we exclude it in future theorems to make their proofs slightly more streamlined.

\begin{thm}\label{th:thicken}{\rm (Stretching in $y$)}
Let~$1<x<y$ and $c \equiv 0 \text{ or } 1 \pmod{y}$.  If~$a+b = y-1$ and~$a \geq \omega(x,b)$ then there is a standard realization for $\{1^a, x^b, y^c\}$. 
\end{thm}

\begin{proof}
By Theorem~\ref{th:omega}, there is a standard linear realization~$\mathbf{h} = [h_1, h_2,  \ldots, h_{a+b+1}]$ for~$\{1^a, x^b\}$.  Let $\boldell = (\ell_1, \ldots, \ell_{a+b})$ be the sequence of differences.
Then 
$$ \apph{\boldell}{k=0}{y-1} \Psi_{h_{k+1}} {\textrm \ \ \ \ and \ \ \ \ } \Psi_0 \ \uplus^{\ell_1} \ 
\apph{\boldell}{k=0}{y-1} \Psi_{y-1-h_{k+1}}$$
give the required realizations for~$c \equiv 0$ and~$1 \pmod{y}$ respectively.

In each case we move through the fauxsets as dictated by the standard linear realization~$\mathbf{h}$, ``forwards" in the first case and ``backwards" in the second.  Note that if we set~$\mathbf{h}$ to be the perfect linear realization~$[0,1,2,\ldots,a+b]$ for~$\{1^{a+b}\}$ we recover the standard linear realizations~$\mathbf{h_1}$ and~$\mathbf{h_2}$.  

If we set up the grid-graph realization so that the tops of all of the fauxsets are on the same level, then the bottoms are also all on the same level, with the exception of~$\varphi_0$ when~$c \equiv 1 \pmod{y}$.  Then all transitions between fauxsets are horizontal (as we do not make a fauxset transition from the vertex~0) and we see that the differences in the construction from moving between fauxsets are indeed exactly those in~$\boldell$.  Internally to the fauxsets we realize~$\{y^c\}$, hence the construction realizes~$\{1^a,x^b,y^c\}$ as required.

The first diagram of Figure~\ref{fig:fullfauxset} illustrates the construction when~$c \equiv 1 \pmod{y}$, using~$x=3$,~$y=7$,~$c=15$ and the standard realization~$[0,3,6,5,2,1,4]$ for~$\{1^2, 3^4\}$.
\end{proof}

The next step is to extend the idea of Theorem~\ref{th:thicken} to other congruence classes of~$v$.  The main issue to deal with is that we now have to bridge between fauxsets of differing sizes, which is not automatically allowable.    

Lemma~\ref{lem:fullfaux} constructs realizations with full fauxsets, based on the $\omega$-constructions~$\mathbf{h_1}$ and~$\mathbf{h_2}$.  Lemma~\ref{lem:tailcurl} gives the construction we use when a full fauxset one is not available, based on the $\omega$-construction~$\mathbf{h_3}$.



\begin{lem}\label{lem:fullfaux}{\rm (Shared-fauxset concatenation)}
For $i \in \{1,2\}$,  let $L_i = \{1^{a_i}, x^{b_i}\}$ be a multiset of size~$v_i = a_i+b_i+1$.  Let~$a = a_1 + a_2$, $b= b_1+b_2$ and $y = a+b+1$.   Suppose~$L_1$ has a perfect linear realization and~$L_2$ has a standard linear realization.

Then there is a standard realization for~$\{1^a, x^b, y^c\}$ of size~$v-1$ in each of the following cases:
\begin{itemize}
\item $v_1$ is even and $v \equiv v_1 \pmod{y}$,
\item $v_1$ is odd and $v \equiv v_2 + 1 \pmod{y}$.
\end{itemize}
\end{lem}

\begin{figure}[tp]
\caption{A standard linear realization for~$\{1^5, 3^{10}, 16^{42} \}$ from the proof of Lemma~\ref{lem:fullfaux}.}\label{fig:1316}
\begin{center}
\begin{tikzpicture}[scale=0.9, every node/.style={transform shape}]
\fill (0,1) circle (2pt) ;\fill (0,2) circle (2pt) ;\fill (0,3) circle (2pt) ;\fill (0,4) circle (2pt) ;\fill (1,1) circle (2pt) ;\fill (1,2) circle (2pt) ;\fill (1,3) circle (2pt) ;\fill (1,4) circle (2pt) ;\fill (2,1) circle (2pt) ;\fill (2,2) circle (2pt) ;\fill (2,3) circle (2pt) ;\fill (2,4) circle (2pt) ;\fill (3,1) circle (2pt) ;\fill (3,2) circle (2pt) ;\fill (3,3) circle (2pt) ;\fill (3,4) circle (2pt) ;\fill (4,1) circle (2pt) ;\fill (4,2) circle (2pt) ;\fill (4,3) circle (2pt) ;\fill (4,4) circle (2pt) ;\fill (5,1) circle (2pt) ;\fill (5,2) circle (2pt) ;\fill (5,3) circle (2pt) ;\fill (5,4) circle (2pt) ;\fill (6,1) circle (2pt) ;\fill (6,2) circle (2pt) ;\fill (6,3) circle (2pt) ;\fill (6,4) circle (2pt) ;\fill (7,1) circle (2pt) ;\fill (7,2) circle (2pt) ;\fill (7,3) circle (2pt) ;\fill (7,4) circle (2pt) ;\fill (8,1) circle (2pt) ;\fill (8,2) circle (2pt) ;\fill (8,3) circle (2pt) ;\fill (8,4) circle (2pt) ;\fill (9,1) circle (2pt) ;\fill (9,2) circle (2pt) ;\fill (9,3) circle (2pt) ;\fill (9,4) circle (2pt) ;\fill (10,1) circle (2pt) ;\fill (10,2) circle (2pt) ;\fill (10,3) circle (2pt) ;\fill (11,1) circle (2pt) ;\fill (11,2) circle (2pt) ;\fill (11,3) circle (2pt) ;\fill (12,1) circle (2pt) ;\fill (12,2) circle (2pt) ;\fill (12,3) circle (2pt) ;\fill (13,1) circle (2pt) ;\fill (13,2) circle (2pt) ;\fill (13,3) circle (2pt) ;\fill (14,1) circle (2pt) ;\fill (14,2) circle (2pt) ;\fill (14,3) circle (2pt) ;\fill (15,1) circle (2pt) ;\fill (15,2) circle (2pt) ;\fill (15,3) circle (2pt) ;
\draw (0,1)--(0,4)--(1,4)--(1,1) ;\draw (2,1)--(2,4)--(3,4)--(3,1) ;\draw (4,1)--(4,4);\draw (5,1)--(5,4);\draw (6,1)--(6,4);\draw  (7,4)--(7,1)--(8,1)--(8,4);\draw (9,1)--(9,4);\draw (10,3)--(10,1)--(11,1)--(11,3);\draw (12,1)--(12,3);\draw (13,1)--(13,3);\draw (14,3)--(14,1)--(15,1)--(15,3);\draw  plot [smooth] coordinates {(1,1) (2.5,0.5) (4,1)};\draw  plot [smooth] coordinates {(2,1) (3.5,0.5) (5,1)};\draw  plot [smooth] coordinates {(3,1) (4.5,0.5) (6,1)};\draw  plot [smooth] coordinates {(9,1) (10.5,0.5) (12,1)};\draw  plot [smooth] coordinates {(4,4) (5.5,4.5) (7,4)};\draw  plot [smooth] coordinates {(5,4) (6.5,4.5) (8,4)};\draw  plot [smooth] coordinates {(6,4) (7.5,4.5) (9,4)};\draw  plot [smooth] coordinates {(10,3) (11.5,3.5) (13,3)};\draw  plot [smooth] coordinates {(11,3) (12.5,3.5) (14,3)}; \draw  plot [smooth] coordinates {(12,3) (13.5,3.5) (15,3)};
\node at (-0.3, 1) {\tiny 0} ; \node at (-0.3, 4) {\tiny 48} ; \node at (9.3, 4) {\tiny 57} ;  \node at (15.3, 3) {\tiny 47} ;  \node at (15.3, 1) {\tiny 15} ;  \node at (13, 0.7) {\tiny 13} ;  
\end{tikzpicture}\end{center}
\end{figure}
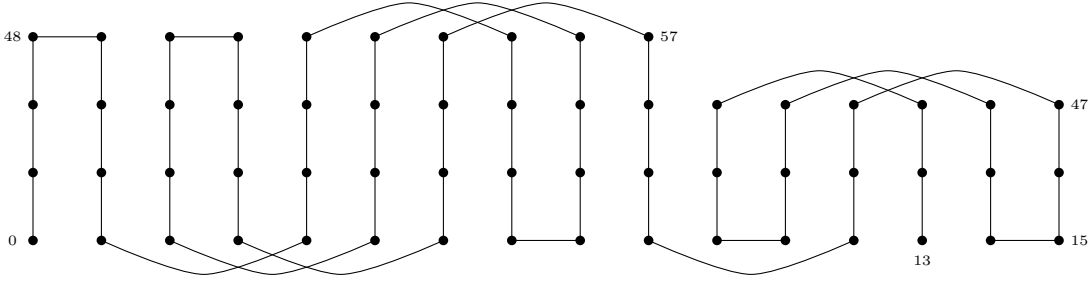

\begin{proof}
The first case is based on~$\mathbf{h_1}$ and the second on~$\mathbf{h_2}$.  In both cases the method is to start by making fauxset hops using the perfect linear realization for~$L_1$, which will use all of one size of fauxset, and follow this with fauxset hops using the standard linear realization for~$L_2$. 

Let~$\mathbf{g_1}$ be the perfect linear realization of~$L_1$ and let~$\mathbf{g_2}$ be the standard linear realization of~$L_2$.
Let~
$$\mathbf{h} = \mathbf{g_1} \oplus \mathbf{g_2} = (h_1, \ldots, h_{y}),$$ 
which is a standard realization for~$\{1^a, x^b\}$.  Denote the 
edge-length sequence $\boldell = (\ell_1, \ldots, \ell_{y-1})$.

Consider the first case, so~$v_1$ is even and~$a+b+c+1 = v \equiv v_1\pmod{y}$.  We claim that the sequence
$$ \apph{\boldell}{k=0}{y-1} \Psi_{h_{k+1}} $$
has the required properties.

The differences are~$\{1^a, x^b, y^c\}$ by construction.  All that needs to be checked is that we have a valid sequence; that is, all the transitions between fauxsets are ``horizontal," in the sense that they connect two lowest elements or two highest ones.  It is sufficient to check that when fauxsets of differing sizes are connected that it happens at the lowest elements. 

The first~$v_1$ fauxsets, which have equal sizes, are traversed using the portion of~$\mathbf{h}$ corresponding to~$\mathbf{g_1}$.  As~$v_1$ is even, this portion of the sequence ends on the lowest element of the fauxset~$\phi_{v_1  -1}$ (which happens to be~$v_1-1$).  The final~$v_2-1$ fauxsets have equal sizes.  Starting from the bottom of~$\phi_{v_1 - 1}$ these are traversed using the portion of~$\mathbf{h}$ corresponding to~$\mathbf{g_2}$. 

The second and third claims follow in the same way, using the sequence
$$ \Psi_0 \ \uplus^{\ell_1} \ 
\apph{\boldell'}{k=0}{y-1} \Psi_{y-1-h_{k+1}},$$
where $\boldell' = (\ell_2, \ldots, \ell_{y-1})$, which works through the fauxsets in decreasing, rather than increasing, order.

Figure~\ref{fig:1316} gives a standard linear realization for~$\{1^5, 3^{10}, 16^{42} \}$ constructed using this method.  We have~$v=58$ and take~$v_1 = 10$, which is even and has~$v \equiv v_1 \pmod{16}$.  We therefore use the first case, based on~$\mathbf{h_1}$. 
 For the perfect linear realization it uses one for~$\{1^3, 3^6\}$ constructed by concatenating the perfect linear realization~$[0,1]$ with the (perfect linear) $\omega$-construction~$\mathbf{h_1}$ for~$\{1^2, 3^6\}$.  For the standard linear realization it uses the $\omega$-construction~$\mathbf{h_2}$ for~$\{1^2, 3^4 \}$.
\end{proof}

\begin{lem}\label{lem:tailcurl}{\rm (Shared-fauxset with tail curl)}
For $i \in \{1,2\}$,  let $L_i = \{1^{a_i}, x^{b_i}\}$ be a multiset of size~$v_i = a_i+b_i+1$.  Let~$a = a_1 + a_2$, $b= b_1+b_2$ and $y = a+b+1$.   Suppose~$L_1$ has a perfect linear realization and~$L_2$ has a standard linear realization whose final edge is a $z$-edge, for $z \in \{1,x\}$.

Then there is a standard realization for~$\{1^a, x^b, y^c\} \cup \{z\}$ of size~$v-1$ in each of the following cases:
\begin{itemize}
\item $v_1$ is even and $v \equiv v_1 \pmod{y}$,
\item $v_1$ is odd and $v \equiv v_2 + 1 \pmod{y}$.
\end{itemize}
\end{lem}

\begin{proof}
The construction is based on the $\omega$-construction~$\mathbf{h_3}$.  If~$v_1$ is odd, as in the proof of Lemma~\ref{lem:fullfaux}, 
let~$\mathbf{g_1}$ be the perfect linear realization of~$L_1$ and let~$\mathbf{g_2}$ be the standard linear realization of~$L_2$.
Let~
$$\mathbf{h} = \mathbf{g_1} \oplus \mathbf{g_2} = (h_1, \ldots, h_{y}),$$ 
which is a standard realization for~$\{1^a, x^b\}$.  Denote the 
edge-length sequence $\boldell = (\ell_1, \ldots, \ell_{y-1})$ and let~$\boldell' = (\ell_2, \ldots, \ell_{y-3})$. 

Consider the sequence
$$\Psi_0 \uplus^{\ell_1} \left( \apph{\boldell'}{k=y-1}{3} \Psi_k \right) \uplus^{\ell_{y-2}} \Psi_2^{(q^*)} \uplus^{\ell_{y-1}} \Psi_1 \uplus^{\ell_{y-1}} \Psi_2^{(0, q^*-1)}. $$ 

As the edge of length~$z = \ell_{y-1}$ is used twice between fauxsets, the differences are $\{1^a, x^b, y^c\} \cup \{z\}$.  We have to check that the sequence is valid.  This follows via an analogous argument to that in the proof of Lemma~\ref{lem:fullfaux}.

Figure~\ref{fig:1316b} gives a standard linear realization for~$\{1^5, 3^{11}, 16^{42} \}$ constructed using this method.  We have~$v=59$ and take~$v_1 = 7$, which is odd, and~$v_2=10$, and so~$v \equiv v_2 +1 \pmod{16}$.  
For the perfect linear realization it uses one for~$\{1^3, 3^3\}$ constructed by concatenating the perfect linear realization~$[0,1]$ with the (perfect linear) $\omega$-construction~$\mathbf{h_1}$ for~$\{1^2, 3^3\}$.  For the standard linear realization it uses the $\omega$-construction~$\mathbf{h_2}$ for~$\{1^2, 3^7 \}$, whose final edge is a 3-edge and so~$z=3$.

If~$v_1$ is even, following the corresponding steps in Lemma~\ref{lem:fullfaux}, with a similar tail curl at the end yields the same result.

A tail curl applied to Figure~\ref{fig:1316} would give a standard linear realization for~$\{1^5, 3^{11}, 16^{41} \}$ by creating a 3-edge between vertices labelled 10 and 13, while removing the 16-edge between vertices labelled 10 and 26.
\end{proof}

\begin{figure}[tp]
\caption{A standard linear realization for~$\{1^5, 3^{11}, 16^{42} \}$ from the proof of Lemma~\ref{lem:tailcurl}.}\label{fig:1316b}
\begin{center}
\begin{tikzpicture}[scale=0.9, every node/.style={transform shape}]
\fill (0,1) circle (2pt) ;\fill (0,2) circle (2pt) ;\fill (0,3) circle (2pt) ;\fill (0,4) circle (2pt) ;\fill (1,1) circle (2pt) ;\fill (1,2) circle (2pt) ;\fill (1,3) circle (2pt) ;\fill (1,4) circle (2pt) ;\fill (2,1) circle (2pt) ;\fill (2,2) circle (2pt) ;\fill (2,3) circle (2pt) ;\fill (2,4) circle (2pt) ;\fill (3,1) circle (2pt) ;\fill (3,2) circle (2pt) ;\fill (3,3) circle (2pt) ;\fill (3,4) circle (2pt) ;\fill (4,1) circle (2pt) ;\fill (4,2) circle (2pt) ;\fill (4,3) circle (2pt) ;\fill (4,4) circle (2pt) ;\fill (5,1) circle (2pt) ;\fill (5,2) circle (2pt) ;\fill (5,3) circle (2pt) ;\fill (5,4) circle (2pt) ;\fill (6,1) circle (2pt) ;\fill (6,2) circle (2pt) ;\fill (6,3) circle (2pt) ;\fill (6,4) circle (2pt) ;\fill (7,1) circle (2pt) ;\fill (7,2) circle (2pt) ;\fill (7,3) circle (2pt) ;\fill (7,4) circle (2pt) ;\fill (8,1) circle (2pt) ;\fill (8,2) circle (2pt) ;\fill (8,3) circle (2pt) ;\fill (8,4) circle (2pt) ;\fill (9,1) circle (2pt) ;\fill (9,2) circle (2pt) ;\fill (9,3) circle (2pt) ;\fill (9,4) circle (2pt) ;\fill (10,1) circle (2pt) ;\fill (10,2) circle (2pt) ;\fill (10,3) circle (2pt) ;\fill (11,1) circle (2pt) ;\fill (11,2) circle (2pt) ;\fill (11,3) circle (2pt) ;\fill (12,1) circle (2pt) ;\fill (12,2) circle (2pt) ;\fill (12,3) circle (2pt) ;\fill (13,1) circle (2pt) ;\fill (13,2) circle (2pt) ;\fill (13,3) circle (2pt) ;\fill (14,1) circle (2pt) ;\fill (14,2) circle (2pt) ;\fill (14,3) circle (2pt) ;\fill (15,0) circle (2pt) ;\fill (15,1) circle (2pt) ;\fill (15,2) circle (2pt) ;\fill (15,3) circle (2pt) ;
\draw (0,4)--(0,1)--(1,1)--(1,4) ;\draw (2,1)--(2,4) ;\draw (3,4)--(3,1) ;\draw (4,1)--(4,4);\draw (5,2)--(5,4);\draw (6,1)--(6,4);\draw  (7,1)--(7,4)--(8,4)--(8,1);\draw (9,1)--(9,4);\draw (10,1)--(10,3)--(11,3)--(11,1);\draw (12,1)--(12,3)--(13,3)--(13,1);\draw (14,1)--(14,3)--(15,3)--(15,0);
\draw  plot [smooth] coordinates {(0,4) (1.5,4.5) (3,4)};\draw  plot [smooth] coordinates {(1,4) (2.5,4.5) (4,4)};\draw  plot [smooth] coordinates {(2,4) (3.5,4.5) (5,4)};\draw  plot [smooth] coordinates {(6,4) (7.5,4.5) (9,4)};\draw  plot [smooth] coordinates {(2,1) (3.5,0.5) (5,1)};\draw  plot [smooth] coordinates {(3,1) (4.5,0.5) (6,1)};\draw  plot [smooth] coordinates {(4,1) (5.5,0.5) (7,1)};\draw  plot [smooth] coordinates {(5,1) (6.5,0.5) (8,1)};\draw  plot [smooth] coordinates {(9,1) (10.5,0.5) (12,1)};\draw  plot [smooth] coordinates {(10,1) (11.5,0.5) (13,1)};\draw  plot [smooth] coordinates {(11,1) (12.5,0.5) (14,1)};
\node at (-0.3, 1) {\tiny 1} ; \node at (-0.3, 4) {\tiny 49} ; \node at (9.3, 4) {\tiny 58} ;  \node at (15.3, 3) {\tiny 48} ;  \node at (15.3, 0) {\tiny 0} ;  \node at (5, 1.7) {\tiny 22} ;  
\end{tikzpicture}\end{center}
\end{figure}

Theorem~\ref{th:odd_hops} applies the methodology to cases with odd~$x$. It complements a similar result for even~$x$ in Theorem~\ref{th:known}.\ref{th:evenx}.

\begin{thm}\label{th:odd_hops}{\rm (Odd $x$)}
Let $x$ be odd and $y > 2x-2$.  The multiset~$\{1^a, x^b, y^c\}$ has a linear realization when $a \geq 4x-3$ and~$b \geq y-2x+2$.
\end{thm}

\begin{proof}
The strategy is to concatenate up to three linear realizations: an $\omega$-construction with support~$\{1,x\}$, a perfect linear realization with support~$\{1\}$, and a standard linear realization with support~$\{1,x,y\}$ constructed from Lemma~\ref{lem:fullfaux} or~\ref{lem:tailcurl}.  The first two of these might be trivial.


Assign~$x$ of the 1-edges to be available to the first realization. (We might only need~$x-1$ or~0, in which cases the second realization can absorb the excess.) 
The first realization can contain as many~$x$-edges as necessary without requiring more $1$-edges.  Similarly, we may continue to add 1-edges to the second realization without constraint.  It therefore suffices to show that~$(4x-3) - x = (3x-3)$ 1-edges and~$(y-2x)$ $x$-edges are sufficient to apply one of Lemma~\ref{lem:fullfaux} or~\ref{lem:tailcurl} to obtain the third realization. 

The total number of vertices for the third realization is $v' = c + \omega(y,c) + 1$; say~$v' = qy+r$ with~$0 \leq r < y$.  
This linear realization, whether constructed from Lemma~\ref{lem:fullfaux} or~\ref{lem:tailcurl}, comes from the concatenation of a perfect linear realization and a standard one, which is then used to guide the path through the fauxsets.  
Let~$v_1$ be the number of vertices available for the first of the two and~$v_2$ be the number available for the second.  So, $\{v_1, v_2\} = \{ r, y-r+1 \}$.  In each case we want to guarantee that there are enough 1-edges and enough $x$-edges to construct the required realization.  

If either $v_1$ or~$v_2$ is less than~$x$, then we can use only 1-edges for that part of the realization, which is more efficient than the general argument of the next two paragraphs.  So assume~$v_1, v_2 \geq x$.

The $\omega$-constructions give perfect realizations with support~$\{1,x\}$ when the number of vertices is congruent to~$0,2 \pmod{x}$ and use~$(x-1)$ 1-edges.  Thus $(x-3) + (x-1) = (2x-4)$ 1-edges and $(v_1 -1) - (x-1) = (v_1-x)$ $x$-edges is always sufficient to construct a perfect linear realization on~$v_1$ vertices. (As discussed above, any excess edges, whether 1-edges or $x$-edges, can be accommodated in the first two realizations.)

It is always possible to construct a standard linear realization on~$v_2$ vertices if we have~$x$ 1-edges and $(v_2-1) - (x-1) = (v_2 - x)$ $x$-edges (with exactly one edge to be accommodated in the first two realizations).  We also need to account for the $z$-edge when using Lemma~\ref{lem:tailcurl}; to guarantee we can always do this we add~1 more of each edge to give $(x+1)$ 1-edges and   $(v_2 - x +1)$ $x$-edges.

Thus the required third realization can always be constructed when we have $(2x-4) + (x+1) = (3x-3)$ 1-edges and $(v_1 - x) + (v_2 - x + 1) = (y - 2x+2)$ $x$-edges.
\end{proof}

Although Theorem~\ref{th:odd_hops} only requires~$y > 2x-2$, given~Theorem~\ref{th:1xy}, it only starts to cover new multisets when~$4x-3 < x+y$.  This applies when $y > 3x-3$.  More generally, the larger~$y$ is compared to~$x$, the more of the target region is covered and so the more useful the result is.  Example~\ref{ex:1532} gives an illustration with~$y \approx 6x$.

\begin{exa}\label{ex:1532}{\rm (Application illustration)}
Let's consider the multiset~$L = \{1^a ,5^b ,32^c\}$ with~$a,b,c \geq 1$.  By Theorem~\ref{th:odd_hops}, we can realize~$L$ whenever $a \geq 4\cdot5 -3 = 17$ and $b \geq 32 - 2\cdot5 = 24$.  This covers many instances not covered by Theorem~\ref{th:1xy}, which says that~$L$ is realizable when~$a \geq 5 + 32 = 37$.

For example, consider~$v = 247$.  The equivalent multisets are $L' = \{ 1^c , 23^b, 54^a \}$ and $L'' = \{ 1^b, 43^c, 99^a \}$.  By Theorem~\ref{th:1xy}, any counterexample must have~$a \leq 36$ (as in the previous paragraph), $b \leq 141$ and $c \leq 76$.  This implies 
$$b  = (v-1) - (a+c) \geq (v-1) - 112 = 134$$
and
$$a = (v-1) - (b+c) \geq (v-1) - 217 = 29.$$
Hence, by Theorem~\ref{th:odd_hops}, no counterexample exists and the Coprime BHR Conjecture holds for the support~$\{1,5,32\}$ when~$v = 247$.

As a second, perhaps more typical, example, consider~$v = 243$.  The equivalent multisets are $L' = \{ 1^c, 38^a, 53^b \}$ and $\{1^b, 55^c, 97^a\}$.  Using the same reasoning as the~$v=247$ case, we find that any counterexample to the Coprime BHR Conjecture must have~$a \leq 36$, $b\leq 151$ and $c \leq 90$.  A quick computer enumeration shows that there are~666 instances satisfying these parameters.  Of these, 530 (about 80\%)  have~$a \geq 17$ and~$b \geq 24$ and so are realizable by Theorem~\ref{th:odd_hops}. This reduces the number of instances for which we do not have a realization to~136. 
\end{exa}

The results presented so far in this section are good general tools for constructing realizations in the target region.   In some situations it is possible to do better still by having edges that hop back and forth between the fauxsets of differing sizes, as in the second realization given in Figure~\ref{fig:fullfauxset}.  It seems that the worst case using this more complicated method does not improve significantly on Theorem~\ref{th:odd_hops} in general, so we do not pursue it in full detail.

However, we conclude this section with some more specialized results and examples where we judge the simplicity to effectiveness ratio to be in our favor.  This gives  Theorem~\ref{th:1xy_new}, which improves on Corollary~\ref{cor:1xy} (and does so to a greater extent than simply applying Theorem~\ref{th:1xy2}), and Theorem~\ref{th:xoddyeven} for the situation when~$x$ is odd and~$y$ is even.

As a strengthened version of Corollary~\ref{cor:1xy}, Theorem~\ref{th:1xy_new} is central to the overarching vision of the paper and will be useful in Sections~\ref{sec:1xtx} and~\ref{sec:conseq}.  Although the number of new realizations is small compared to other results in the paper, they are guaranteed to hit parameters in the target region (unless the target region is already empty).

\begin{thm}\label{th:1xy_new}{\rm (Counterexample characterization refined)}
Let $L = \{1^a, x^b, y^c\}$ be strongly admissible with equivalent multisets $L' = \{1^b, \dot{x}^c, \dot{y}^a  \}$ and $L'' = \{1^c, \ddot{x}^a, \ddot{y}^b\}$.   Define a function
$$f(x,y) = \begin{cases}
  y-1 & \text{if $x$ and $y$ are even} \\
  x+y-2 & \text{if $x=3$ or if $x$ is even and $y$ is odd,} \\
  x+y-1 & \text{otherwise.}
\end{cases}$$
If~$a\geq f(x,y)$, $b \geq f(\dot{x}, \dot{y})$ or $c \geq f( \ddot{x} , \ddot{y})$ then~$L$ is realizable.  In particular, any counterexample to the Coprime BHR Conjecture has
$$ f(x , y) + f(\dot{x} ,\dot{y}) + f(\ddot{x} , \ddot{y}) \geq v+2. $$
\end{thm}

\begin{proof}
This almost follows from Theorem~\ref{th:1xy2} in the same way that Corollary~\ref{cor:1xy} follows from Theorem~\ref{th:1xy}. 
Additional work is needed when~$y$ is odd.  In particular,
when~$x=3$ then we need to lower the bound on~$a$ by 2 and when $x$ is odd with~$x>3$ we need to lower it by~1. 

The first additional piece needed is the observation that a multiset with support~$\{1,3\}$ always has a standard linear realization that uses two 1-edges via Theorem~\ref{th:omega}.  This gives a bound of~$a \geq x+y-1$  when~$x=3$ and~$y$ is odd, leaving an improvement of~1 in the bound needed in all cases with~$x$ odd.

The subcases that require this improvement are those for which a tail curl is needed and the~$\omega$-construction has the form~$\{1^y, y^c\}$  To complete the proof, we construct a standard linear realization for $\{1^{y-1}, x, y^c\}$ for these subcases, which are then concatenated with a standard linear realization for~$\{1^{a-y+1}, x^{b-1}\}$ to give the result.  (Note that we may assume that~$b \geq 1$ as otherwise we do not need a concatenation and have the stronger bound~$a \geq y$.)

Writing~$c = q'y+r'$ with $0 \leq r' < y$, we have that each of~$x$,~$y$ and~$r'$ are odd and~$q' > 0$ in the problematic subcases (as are stronger conditions but these are sufficient for the construction to succeed). 

There are two constructions, one for $x > y-r'$ and one for~$x < y-r'$ (note that $y-r'$ is even, so we never have $x = y-r'$).

First, suppose that~$x > y-r'$. Let
$$\mathbf{g_1}  = \left(  \apph{}{k=0}{y-x-1} \Psi_{k} \right)
\ \uplus^{x} \ 
\left( \apph{}{k=y-1}{y-x} \Psi_{k} \right),$$ 
with first element~0.

We claim that~$\mathbf{g_1}$ is a standard linear realization for $\{1^{y-2}, x, y^c\}$.  The differences are correct by construction, provided it is a valid sequence.  There are two transitions between fauxsets of differing sizes to check.  

First, the $x$-edge from fauxset~$\phi_{y-x-1}$ to~$\phi_{y-1}$.  As both~$x$ and~$y$ are odd, $y-x-1$ is also odd.  Therefore we leave~$\phi_{y-x-1}$ from its lowest element (which is~$y-x-1$) and we may connect with an $x$-edge to the lowest element of~$\phi_{y-1}$ (which is~$y-1$).

Second, the $1$-edge from fauxset~$\phi_{r'}$ to fauxset~$\phi_{r'-1}$, which happens within the second part of~$\mathbf{g_1}$ after the $x$-edge has been used.  As~$r'$ is odd, we leave~$\phi_{r'}$ from its lowest element (which is~$r'$) and we may connect with an $1$-edge to the lowest element of~$\phi_{r'-1}$ (which is~$r'-1$).

The concatenation~$[0,1] \oplus \mathbf{g_1}$ is then the required realization for $\{1^{y-1}, x, y^c\}$.
The first diagram in Figure~\ref{fig:1xedge} illustrates~$\mathbf{g_1}$ for the multiset~$\{1^7,5,9^{25}\}$.  

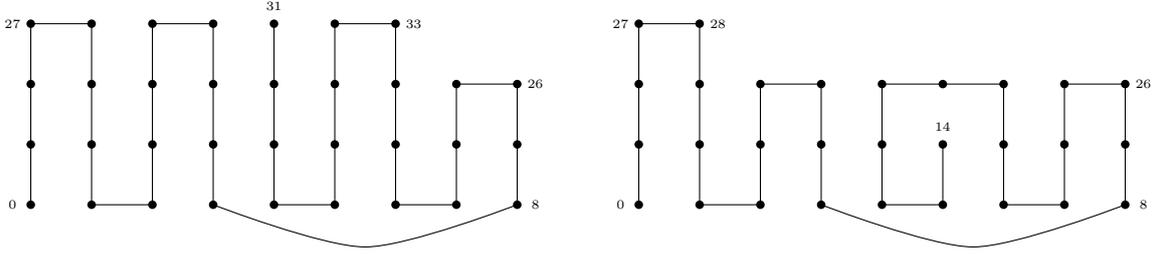
\begin{figure}[tp]
\caption{Standard linear realizations for~$\{1^7, 5, 9^{25} \}$ and~$\{1^8, 5, 9^{19} \}$ from the proof of Theorem~\ref{th:1xy_new}.}\label{fig:1xedge}
\begin{center}
\begin{tikzpicture}[scale=0.8, every node/.style={transform shape}]
\fill (0,1) circle (2pt) ;\fill (0,2) circle (2pt) ;\fill (0,3) circle (2pt) ;\fill (0,4) circle (2pt) ;\fill (1,1) circle (2pt) ;\fill (1,2) circle (2pt) ;\fill (1,3) circle (2pt) ;\fill (1,4) circle (2pt) ;\fill (2,1) circle (2pt) ;\fill (2,2) circle (2pt) ;\fill (2,3) circle (2pt) ;\fill (2,4) circle (2pt) ;\fill (3,1) circle (2pt) ;\fill (3,2) circle (2pt) ;\fill (3,3) circle (2pt) ;\fill (3,4) circle (2pt) ;\fill (4,1) circle (2pt) ;\fill (4,2) circle (2pt) ;\fill (4,3) circle (2pt) ;\fill (4,4) circle (2pt) ;\fill (5,1) circle (2pt) ;\fill (5,2) circle (2pt) ;\fill (5,3) circle (2pt) ;\fill (5,4) circle (2pt) ;\fill (6,1) circle (2pt) ;\fill (6,2) circle (2pt) ;\fill (6,3) circle (2pt) ;\fill (6,4) circle (2pt) ;\fill (7,1) circle (2pt) ;\fill (7,2) circle (2pt) ;\fill (7,3) circle (2pt) ;\fill (8,1) circle (2pt) ;\fill (8,2) circle (2pt) ;\fill (8,3) circle (2pt) ;\fill (10,1) circle (2pt) ;\fill (10,2) circle (2pt) ;\fill (10,3) circle (2pt) ;\fill (10,4) circle (2pt) ;\fill (11,1) circle (2pt) ;\fill (11,2) circle (2pt) ;\fill (11,3) circle (2pt) ;\fill (11,4) circle (2pt) ;\fill (12,1) circle (2pt) ;\fill (12,2) circle (2pt) ;\fill (12,3) circle (2pt) ;\fill (13,1) circle (2pt) ;\fill (13,2) circle (2pt) ;\fill (13,3) circle (2pt) ;\fill (14,1) circle (2pt) ;\fill (14,2) circle (2pt) ;\fill (14,3) circle (2pt) ;\fill (15,1) circle (2pt) ;\fill (15,2) circle (2pt) ;\fill (15,3) circle (2pt) ;\fill (16,1) circle (2pt) ;\fill (16,2) circle (2pt) ;\fill (16,3) circle (2pt) ;\fill (17,1) circle (2pt) ;\fill (17,2) circle (2pt) ;\fill (17,3) circle (2pt) ;\fill (18,1) circle (2pt) ;\fill (18,2) circle (2pt) ;\fill (18,3) circle (2pt) ;
\draw (0,1)--(0,4)--(1,4)--(1,1)--(2,1)--(2,4)--(3,4)--(3,1) ;\draw (10,1)--(10,4)--(11,4)--(11,1)--(12,1)--(12,3)--(13,3)--(13,1) ;\draw (8,1)--(8,3)--(7,3)--(7,1)--(6,1)--(6,4)--(5,4)--(5,1)--(4,1)--(4,4);\draw (18,1)--(18,3)--(17,3)--(17,1)--(16,1)--(16,3)--(14,3)--(14,1)--(15,1)--(15,2) ;
\draw  plot [smooth] coordinates {(3,1) (5.5,0.3) (8,1)}; \draw  plot [smooth] coordinates {(13,1) (15.5,0.3) (18,1)};
\node at (-0.3, 1) {\tiny 0} ; \node at (-0.3, 4) {\tiny 27} ; \node at (8.3, 1) {\tiny 8} ;  \node at (8.3, 3) {\tiny 26} ;  \node at (6.3, 4) {\tiny 33} ;  \node at (4, 4.3) {\tiny 31} ;  \node at (9.7, 1) {\tiny 0} ; \node at (9.7, 4) {\tiny 27} ; \node at (18.3, 1) {\tiny 8} ;  \node at (18.3, 3) {\tiny 26} ;  \node at (11.3, 4) {\tiny 28} ;  \node at (15, 2.3) {\tiny 14} ;  
\end{tikzpicture}\end{center}
\end{figure}

Second, suppose that~$x < y-r'$.  Let
$$\mathbf{g_2}  = \left(  \apph{}{k=0}{y-x-1} \Psi_{k} \right)
\ \uplus^{x} \ 
\left( \apph{}{k=y-1}{y-x+2} \Psi_{k} \right) \uplus \Psi_{y-x+1}^{(q^*)} \uplus \Psi_{y-x} \uplus \Psi_{y-x+1}^{(0, q^*-1)}, $$ 
with first element~0.

We claim that~$\mathbf{g_2}$ is a standard linear realization for~$\{1^{y-1}, x, y^c\}$.  Note that, unlike~$\mathbf{g_1}$, we have~$v = y+c$ and so do not want a full fauxset construction.  The idea is to use a tail curl, in a similar manner to~$\mathbf{h_3}$.   With this observation, the differences are correct by construction, if the sequence is valid.  We again have two transitions to check.

This time, the 1-edge transition between fauxsets of differing sizes happens in the first part of the realization, before the $x$-edge; specifically from fauxset~$\phi_r'$ to~$\phi_{r'+1}$.  As~$r'$ is odd, this transition occurs between the lowest elements of these two fauxsets (elements~$r'$ and~$r'+1$).
The $x$-edge transition, and its justification, work exactly as in~$\mathbf{g_1}$.

The second diagram in Figure~\ref{fig:1xedge} illustrates this for the multiset~$\{1^8,5,9^{21}\}$.
\end{proof}

Earlier results in this section are at their best when~$x$ is small.  The next result is useful when~$x$ is closer to~$y$ in size.  In light of the strong results in Theorems~\ref{th:known}.5 and~\ref{th:ao_largev}.2, the exclusion of~$x = y-1$ from this result is not important in practice.

\begin{thm}\label{th:xoddyeven}{\rm (Odd $x$, even $y$)}
Let~$y$ be even and let~$x$ be odd with~$y/2 \leq x < y-1$.  Let~$d$ be even with~$d \leq y-x-1$.  There is a standard linear realization for $\{1^{y-1-d}, x^d, y^c\}$. 
Hence there is a linear realization for~$\{1^a, x^b, y^c\}$ whenever~$a \geq 2x$ and $b \geq y-x-1$. 
\end{thm}

\begin{proof}
The method is to use full fauxsets with respect to~$y$ with 1-edges and $x$-edges between them.  

Let~$\mathbf{g_1}$ be the perfect linear realization $[0,1,\ldots, y-x-d]$ for~$\{1^{y-x-d}\}$ and let~$\mathbf{g_2}$ be the $\mathbf{h_1}$ pattern $\omega$-construction for~$\{1^{x-1}, x^d\}$.  Consider~$\mathbf{g_1} \oplus \mathbf{g_2}$, with difference sequence~$\boldell$.  As~$\mathbf{g_2}$ has fewer than~$2x$ vertices, there are no consecutive $x$-edges.  Furthermore, all~$x$-edges are in odd positions in the difference sequence for~$\mathbf{g_2}$ and so, as~$y-x-d$ is odd, all $x$-edges in $\boldell$ are in even positions.  This information is sufficient to verify that
$$\apph{\boldell}{k=0}{y-1} \Psi_{k}  {\text \ \ \ \ and \ \ \ \ }    \Psi_0 \uplus \left( \apph{\boldell}{k=y-1}{1} \Psi_k \right)$$
are standard linear realizations for $\{1^{y-1-d}, x^d, y^c\}$ when~$c$ is even or odd respectively.

Figure~\ref{fig:grid_eg} and Figure~\ref{fig:11116} show the standard linear realizations constructed with this method for the multisets $\{1^5,5^2,8^{21}\}$ and $\{1^{11}, 11^{4}, 16^{36} \}$ respectively.

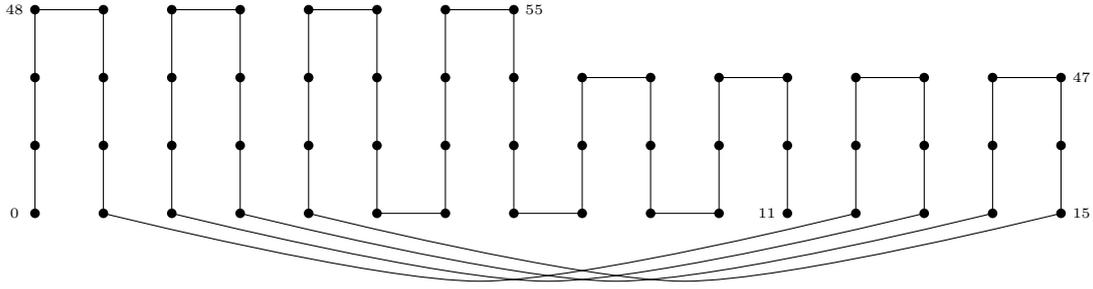
\begin{figure}[tp]
\caption{A standard linear realization for~$\{1^{11}, 11^{4}, 16^{36} \}$ from the proof of Theorem~\ref{th:xoddyeven}.}\label{fig:11116}
\begin{center}
\begin{tikzpicture}[scale=0.9, every node/.style={transform shape}]
\fill (0,1) circle (2pt) ; \fill (0,2) circle (2pt) ; \fill (0,3) circle (2pt) ; \fill (0,4) circle (2pt) ; \fill (1,1) circle (2pt) ; \fill (1,2) circle (2pt) ; \fill (1,3) circle (2pt) ; \fill (1,4) circle (2pt) ; \fill (2,1) circle (2pt) ; \fill (2,2) circle (2pt) ; \fill (2,3) circle (2pt) ; \fill (2,4) circle (2pt) ; \fill (3,1) circle (2pt) ; \fill (3,2) circle (2pt) ; \fill (3,3) circle (2pt) ; \fill (3,4) circle (2pt) ; \fill (4,1) circle (2pt) ; \fill (4,2) circle (2pt) ; \fill (4,3) circle (2pt) ; \fill (4,4) circle (2pt) ; \fill (5,1) circle (2pt) ; \fill (5,2) circle (2pt) ; \fill (5,3) circle (2pt) ; \fill (5,4) circle (2pt) ; \fill (6,1) circle (2pt) ; \fill (6,2) circle (2pt) ; \fill (6,3) circle (2pt) ; \fill (6,4) circle (2pt) ; \fill (7,1) circle (2pt) ; \fill (7,2) circle (2pt) ; \fill (7,3) circle (2pt) ; \fill (7,4) circle (2pt) ; \fill (8,1) circle (2pt) ; \fill (8,2) circle (2pt) ; \fill (8,3) circle (2pt) ; \fill (9,1) circle (2pt) ; \fill (9,2) circle (2pt) ; \fill (9,3) circle (2pt) ; \fill (10,1) circle (2pt) ; \fill (10,2) circle (2pt) ; \fill (10,3) circle (2pt) ; \fill (11,1) circle (2pt) ; \fill (11,2) circle (2pt) ; \fill (11,3) circle (2pt) ; \fill (12,1) circle (2pt) ; \fill (12,2) circle (2pt) ; \fill (12,3) circle (2pt) ; \fill (13,1) circle (2pt) ; \fill (13,2) circle (2pt) ; \fill (13,3) circle (2pt) ; \fill (14,1) circle (2pt) ; \fill (14,2) circle (2pt) ; \fill (14,3) circle (2pt) ; \fill (15,1) circle (2pt) ; \fill (15,2) circle (2pt) ; \fill (15,3) circle (2pt) ;
\draw (0,1)--(0,4)--(1,4)--(1,1) ;\draw (2,1)--(2,4)--(3,4)--(3,1) ;\draw (4,1)--(4,4)--(5,4)--(5,1)--(6,1)--(6,4)--(7,4)--(7,1)--(8,1)--(8,3)--(9,3)--(9,1)--(10,1)--(10,3)--(11,3)--(11,1); \draw (12,1)--(12,3)--(13,3)--(13,1); \draw (14,1)--(14,3)--(15,3)--(15,1);
\draw  plot [smooth] coordinates {(1,1) (6.5,0) (12,1)}; \draw  plot [smooth] coordinates {(2,1) (7.5,0) (13,1)};\draw  plot [smooth] coordinates {(3,1) (8.5,0) (14,1)}; \draw  plot [smooth] coordinates {(4,1) (9.5,0) (15,1)};
\node at (-0.3, 1) {\tiny 0} ; \node at (-0.3, 4) {\tiny 48} ; \node at (7.3, 4) {\tiny 55} ;  \node at (15.3, 3) {\tiny 47} ;  \node at (15.3, 1) {\tiny 15} ;  \node at (10.7, 1) {\tiny 11} ;  
\end{tikzpicture}\end{center}
\end{figure}

To construct a linear realization for $\{1^a, x^b, y^c\}$ with~$a \geq 2x$ and $b \geq y-x-1$, take the above standard linear realization with~$d = y-x-1$, which realizes $\{1^x, x^{y-x-1}, y^c\}$, and concatenate it with a standard linear realization for~$\{1^{a-x} , x^{b - (y-x-1)} \}$, which exists as $$a-x \geq x \geq \omega(x,b-(y-x-1)),$$ to obtain a linear realization for~$\{1^a, x^b, y^c\}$.
\end{proof}

\section{Replacing $x$-edges with $y$-edges}\label{sec:1xtx}

In this section we consider multisets with support~$\{1,x,y\}$ with~$y=tx$.  The main method is to start with an $\omega$-construction for a multiset with support~$\{1,x\}$ and replace some of the~$x$-edges with~$y$-edges.  This results in constructions that apply when the number of~$1$-edges is at least~$x-1$ or~$x$, with the trade-off for this strong low bound being a requirement of a minimum number of~$x$ edges.  

Following these constructions, we give a separate approach based on modular arithmetic considerations that applies for all~$x$. Theorem~\ref{th:1xtx} shows for any fixed support $\{1,x,tx\}$ that multisets with that support are realizable for all sufficiently large~$v$ with~$v \equiv \pm 1 \pmod{tx}$.  When both~$t$ and~$x$ are at least~7, all~$v$ are sufficiently large.

The crucial observation that the constructive method relies on is that when~$y$ is a multiple of~$x$ then a $y$-edge in a linear realization connects vertices in the same fauxset with respect to~$x$.   We can thus think of the fauxset traversals independently of each other and replace $x$-edges with $y$-edges in a given fauxset without reference to the rest of the realization.  We can do better still if we sometimes consider the fauxsets in pairs (Lemma~\ref{lem:pairstd}).

We limit the constructive results to odd~$x$.  As in the last section, existing results for even~$x$ tend to be stronger than are achievable with the new methodology.  In fact, it is even more clear cut here: as $tx$ is even when~$x$ is even, we may use the very strong Theorem~\ref{th:1xy2}.4.

The main result of the section is Theorem~\ref{th:1xtx_gen}, which says, approximately, that $\{1^a, x^b, tx^c\}$ has a linear realization for odd~$x$ provided that $a \geq \omega(x,b+c)$ and $b \geq 2 \lfloor v/x \rfloor + \epsilon$, where $\epsilon$ is a constant that depends only on~$t$ and~$x$.
This requires  two construction tools, given in Lemmas~\ref{lem:endstd} and~\ref{lem:pairstd}. 

This is refined in Theorem~\ref{th:1xtx_tight} to give stronger results, dependent on the parity of~$t$ and the congruence class of~$v$ modulo~$x$.  This needs two further tools, given in Lemmas~\ref{lem:midperf} and~\ref{lem:even_perf}.

\begin{lem}\label{lem:endstd}{\rm (Reordering the terminal fauxset)}
Suppose a standard realization for a multiset~$L$ finishes with the subsequence~$\Psi_k^{(s_1,s_2)}$, a fauxset segment with respect to~$x$ in order with~$w = s_2 - s_1$ edges.  Let~$M$ be a multiset of size~$w$ that has a standard realization.  Then~$(L \setminus \{x^w\}) \cup xM$ has a standard realization.
\end{lem}

\begin{proof}
Suppose the fauxset is traversed in increasing order and
let $$\Psi_k^{(s_1,s_2)} = [s_1x+k, (s_1+1)x + k, \ldots s_2x+k],$$ which realizes~$\{x^w\}$. Let the standard realization for~$M$ be $\mathbf{h} = [h_1, h_2, \ldots, h_{w+1}]$.  Then
$$(\mathbf{h}+s_1)x + k = [(h_1+s_1)x + k, (h_2+s_1)x + k, \ldots, (h_{w+1}x+s_1) + k]$$
realizes~$xM$ and uses the same vertices as~$\Psi_k^{(s_1,s_2)}$.

Replace~$\Psi_k^{(s_1,s_2)}$ with~$(\mathbf{h}+s_1)x+k$ to obtain the desired realization.  A similar argument applies if the original fauxset is traversed in decreasing order.
\end{proof}

\begin{lem}\label{lem:pairstd}{\rm (Reordering a pair of fauxsets)}
Suppose a standard realization for a multiset~$L$ includes a subsequence of the form $\Psi_k^{(s,q^*)} \uplus \Psi_{k+1}^{(s,q^*)}$, where~$\Psi_k$ and~$\Psi_{k+1}$ are fauxsets with respect to~$x$ of the same size and the bridge is at the highest elements, and let~$w=q^*-s$. 
Let~$M$ be a multiset of size~$w$ that has a standard realization.  Then~$(L \setminus \{x^{2w}\}) \cup (xM \cup xM)$ has a standard realization.
\end{lem}

\begin{proof}
Let
$$\Psi_k^{(s,q^*)} \uplus \Psi_{k+1}^{(s,q^*)} = [sx+k, 
\ldots q^*x+k, q^*x + (k+1), 
\ldots, sx+ (k+1)           ],$$ 
which realizes~$\{1, x^{2w}\}$.  
Let the standard realization for~$M$ be $\mathbf{h} = [h_1, h_2, \ldots, h_{s+1}]$.

Replace~$\Psi_k^{(s,q^*)}$ with~$(\mathbf{h}+s)x+k$ and replace~$\Psi_{k+1}^{(s,q^*)}$ with~$(\mathbf{h}+s)x+(k+1)$.  This portion of the realization now realises~$xM \cup \{1\} \cup xM$.  Therefore, the new realization as a whole realizes
$$(L \setminus \{1,x^{2w}\}) \cup (xM \cup \{1\} \cup xM) = (L \setminus \{x^{2w}\}) \cup (xM \cup xM).$$
\end{proof}

\begin{figure}[tp]
\caption{Standard linear realizations for~$\{1^6, 7^{8}, 14^{10}   \}$ and for~$\{1^9, 9^{14}, 18^9 \}$.}\label{fig:1x2x_eg}
\begin{center}
\begin{tikzpicture}[scale=0.8, every node/.style={transform shape}]
\fill (0,1) circle (2pt) ; \fill (0,2) circle (2pt) ; \fill (0,3) circle (2pt) ; \fill (0,4) circle (2pt) ; \fill (1,1) circle (2pt) ; \fill (1,2) circle (2pt) ; \fill (1,3) circle (2pt) ; \fill (1,4) circle (2pt) ; \fill (2,1) circle (2pt) ; \fill (2,2) circle (2pt) ; \fill (2,3) circle (2pt) ; \fill (2,4) circle (2pt) ; \fill (3,1) circle (2pt) ; \fill (3,2) circle (2pt) ; \fill (3,3) circle (2pt) ; \fill (4,1) circle (2pt) ; \fill (4,2) circle (2pt) ; \fill (4,3) circle (2pt) ;\fill (5,1) circle (2pt) ; \fill (5,2) circle (2pt) ; \fill (5,3) circle (2pt) ; \fill (6,0) circle (2pt) ; \fill (6,1) circle (2pt) ; \fill (6,2) circle (2pt) ; \fill (6,3) circle (2pt) ;
\draw (0,3)--(0,4) ; \draw (0,1)--(1,1) ; \draw (1,3)--(1,4) ; \draw (2,3)--(2,4) ; \draw (3,3)--(3,2)--(4,2)--(4,3) ; ; \draw (5,3)--(5,2)--(6,2)--(6,3) ;   \draw (2,1)--(3,1) ; \draw (4,1)--(5,1) ; \draw (6,1)--(6,0) ; \draw (1,2)--(2,2) ; 
\draw  plot [smooth] coordinates {(0,1) (-0.35,2) (0,3)}; \draw  plot [smooth] coordinates {(0,2) (-0.35,3) (0,4)}; \draw  plot [smooth] coordinates {(1,1) (0.65,2) (1,3)}; \draw  plot [smooth] coordinates {(1,2) (0.65,3) (1,4)}; \draw  plot [smooth] coordinates {(2,1) (2.35,2) (2,3)}; \draw  plot [smooth] coordinates {(2,2) (2.35,3) (2,4)}; \draw  plot [smooth] coordinates {(3,1) (2.65,2) (3,3)}; \draw  plot [smooth] coordinates {(4,1) (4.35,2) (4,3)}; \draw  plot [smooth] coordinates {(5,1) (4.65,2) (5,3)}; \draw  plot [smooth] coordinates {(6,1) (6.35,2) (6,3)};
\node at (-0.3, 1) {\tiny 1} ; \node at (0.3, 2) {\tiny 8} ;  \node at (-0.3, 4) {\tiny 22} ;  \node at (6.3, 0) {\tiny 0} ;    \node at (6.3, 1) {\tiny 7} ;  \node at (6.3, 3) {\tiny 21} ;  \node at (2.3, 4) {\tiny 24} ;  
\fill (8,1) circle (2pt) ; \fill (8,2) circle (2pt) ; \fill (8,3) circle (2pt) ; \fill (8,4) circle (2pt) ; \fill (9,1) circle (2pt) ; \fill (9,2) circle (2pt) ; \fill (9,3) circle (2pt) ; \fill (9,4) circle (2pt) ; \fill (10,1) circle (2pt) ; \fill (10,2) circle (2pt) ; \fill (10,3) circle (2pt) ; \fill (10,4) circle (2pt) ; \fill (11,1) circle (2pt) ; \fill (11,2) circle (2pt) ; \fill (11,3) circle (2pt) ; \fill (11,4) circle (2pt) ; \fill (12,1) circle (2pt) ; \fill (12,2) circle (2pt) ; \fill (12,3) circle (2pt) ; \fill (12,4) circle (2pt) ; \fill (13,1) circle (2pt) ; \fill (13,2) circle (2pt) ; \fill (13,3) circle (2pt) ; \fill (14,1) circle (2pt) ; \fill (14,2) circle (2pt) ; \fill (14,3) circle (2pt) ; \fill (15,1) circle (2pt) ; \fill (15,2) circle (2pt) ; \fill (15,3) circle (2pt) ; \fill (16,0) circle (2pt) ; \fill (16,1) circle (2pt) ; \fill (16,2) circle (2pt) ; \fill (16,3) circle (2pt) ;
\draw (9,2) -- (9,3) ;  \draw (11,3) -- (11,4) ; \draw (12,3) -- (12,4) ; \draw (9,1) -- (8,1)   -- (8,4) -- (10,4) -- (10,1) -- (11,1) ;  \draw (11,2) -- (12,2) ; \draw (12,1) -- (13,1) ; \draw  (13,3) -- (13,2) -- (14,2) -- (14,3) ; \draw (14,1) -- (15,1) ; \draw (15,3) --(15,2) -- (16,2) -- (16,3) ; \draw (16,1) -- (16,0) ; \draw  plot [smooth] coordinates {(9,1) (8.65,2) (9,3)};\draw  plot [smooth] coordinates {(11,1) (10.65,2) (11,3)};\draw  plot [smooth] coordinates {(11,2) (10.65,3) (11,4)};\draw  plot [smooth] coordinates {(12,1) (12.35,2) (12,3)};\draw  plot [smooth] coordinates {(12,2) (12.35,3) (12,4)};\draw  plot [smooth] coordinates {(13,1) (12.65,2) (13,3)};\draw  plot [smooth] coordinates {(14,1) (14.35,2) (14,3)};\draw  plot [smooth] coordinates {(15,1) (14.65,2) (15,3)};\draw  plot [smooth] coordinates {(16,1) (16.35,2) (16,3)};
\node at (7.7, 1) {\tiny 1} ;  \node at (7.7, 2) {\tiny 10} ;  \node at (7.7, 4) {\tiny 28} ;  \node at (9.3, 2) {\tiny 11} ;\node at (12.3, 4) {\tiny 32} ;\node at (16.3, 0) {\tiny 0} ;  \node at (16.3, 1) {\tiny 9} ;  \node at (16.3, 3) {\tiny 27} ;  
\end{tikzpicture}\end{center}
\end{figure}
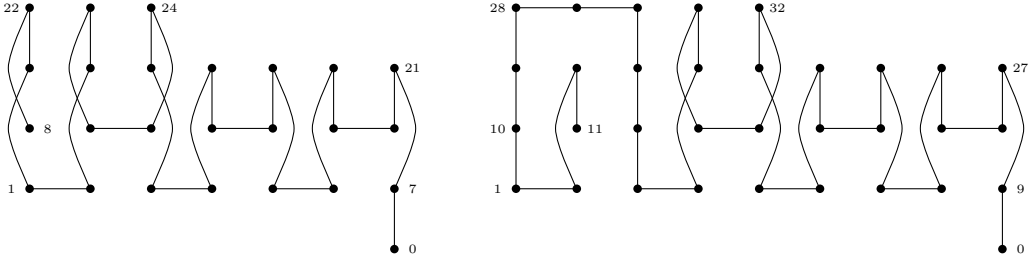

The two realizations in Figure~\ref{fig:1x2x_eg} each show the application of Lemma~\ref{lem:pairstd} to three pairs of adjacent fauxsets and Lemma~\ref{lem:endstd} to the final fauxset to be traversed.

We are now in position to prove the first main result of the section.

\begin{thm}\label{th:1xtx_gen}{\rm (Odd $x$)}
Let~$L = \{1^a, x^b, tx^c \}$ be a multiset of size~$v-1$ with~$x$ odd and $v \geq (t+2)x$.  Let~$q = \lfloor v/x \rfloor$.    If~$a \geq \omega(x,b+c)$ and $b \geq  2q + tx -2t + 1$,
then~$L$ has a standard linear realization.
\end{thm}

\begin{proof}
Start with the $\omega$-construction for~$\{1^{\omega(x,b+c)}, x^{b+c} \}$, with the appropriate constraint on~$b$.  The method is to replace~$c$ of the $x$-edges with $(tx)$-edges.  We accomplish this by showing that any number of $(tx)$-edges can be accommodated provided that the constraints on~$a$ and~$b$ are met.

As~$v \geq (t+2)x$ we have $t < q-1$.  The final fauxset traversal has at least $q-1 > t$ vertices (it can only be as low as~$q-1$ in cases that require a tail curl and that final traversal omits the largest fauxset vertex).  

In the $\omega$-construction~$\mathbf{h_1}$ with support~$\{1,x\}$ it is possible to partition all fauxsets except the final one into adjacent pairs of the same height.   The same is almost true for the $\omega$-construction~$\mathbf{h_2}$: for the construction we leave the edge between~0 and~$x$ as an~$x$-edge and consider the fauxset traversal to start at~$x$---this leads to the ``+1" in the bound calculation below.  Given a pair of such fauxsets with~$q^*$ edges each, we may use Lemma~\ref{lem:pairstd} and a standard realization with support~$\{1,t\}$ to replace any even number~$2m$ of~$x$ edges with~$(tx)$-edges with~$m \leq \omega(t,q^*) \leq t$.   Looked at from the perspective of~$x$-edges, provided we have at least~$t$ $x$-edges per fauxset for these fauxsets and the total number of such~$x$-edges is even, we may complete the realization with~$tx$-edges.

The final fauxset to be traversed uses at least~$t+1$ vertices in this traversal.  We may use Lemma~\ref{lem:endstd} to add at as many~$(tx)$-edges as we like within this fauxset, provided we have at least~$t$ $x$-edges.  Note that in the case with exactly~$t+1$ vertices we only require~$t-1$ $x$-edges to add~1 $(tx)$-edge, so in every case the number of $(tx)$-edges we may add in this fauxset is strictly positive.  This is important as it allows us to control the overall parity of the number of added $(tx)$-edges.

This leaves two fauxsets untouched, with up to $q$ $x$-edges each. 
Therefore, we may add any number of $(tx)$-edges provided that
\begin{eqnarray*}
b & \geq & (x-3)t + 1 + t + 2q \\
  & = & 2q + tx -2t + 1,
\end{eqnarray*}
as required.
\end{proof}

The strength of the results in Theorem~\ref{th:1xtx_gen} might be obscured by number of parameters involved and the slightly complicated dependencies.  However, since~$\omega(x,b+c) \in \{x-1, x\}$, the bound on~$a$ is very strong compared to most known results.  Moreover, for a fixed support~$\{1,x,tx\}$, the proportion~$b/(v-1)$ of~$x$-edges required to meet the constraints tends to~$2/x$  as~$v$ increases.

The standard linear realizations for~$\{1^6, 7^{8}, 14^{10}   \}$ and for~$\{1^9, 9^{14}, 18^9 \}$ in  Figure~\ref{fig:1x2x_eg} are built from the $\omega$-constructions for $\{1^6,7^{18}\}$ and~$\{1^9,9^{23}\}$ (see Figure~\ref{fig:omega_eg}) using the methods in the proof of Theorem~\ref{th:1xtx_gen}. 

These realizations suggest several possibilities for strengthening Theorem~\ref{th:1xtx_gen} further:
\begin{itemize}
\item Other than the final one, there are no unpaired fauxsets necessary in the first realization.  This occurs when no tail curl is required; in this situation we may remove the~$q$-component of the bound.
\item In this example the ``$\omega$ bound" for the subsidiary realizations with support~$\{1,t\}$ is~$t-1$ rather than~$t$ in all fauxsets, meaning that we require one fewer~$x$-edge per fauxset than the general argument.
\item When we have an unused (and unpaired) fauxset we may replace yet more $x$-edges with $(tx)$-edges, but a standard linear realization is not sufficient to do so.  It must be perfect.
\end{itemize}
We pursue the first and third of these observations in Theorem~\ref{th:1xtx_tight} to give improved, albeit more intricate, bounds that guarantee the existence of linear realizations.  (The gains from the second are small, but could still potentially be beneficial in particular cases.)

Two further tools are needed: a perfect linear realization to use when~$t$ is even and a method allowing substitutions in unpaired fauxsets.  These are given in Lemmas~\ref{lem:even_perf} and~\ref{lem:midperf}.

\begin{lem}\label{lem:even_perf}{\rm (Perfect realizations)}
Let~$x$ be even.  There are perfect realizations for $\{ 1^{x+2s-2}, x^{2sx-2s}\}$ and $\{ 1^{x+2s-1} , x^{2sx - 2s + 2} \}$ for each~$s \geq 1$.
\end{lem}

\begin{figure}[tp]
\caption{Perfect linear realizations for~$\{1^7, 6^{12}  \}$ and for~$\{1^8, 6^{20} \}$ from the proof of Lemma~\ref{lem:even_perf}.}\label{fig:even_perf}
\begin{center}
\begin{tikzpicture}[scale=0.9, every node/.style={transform shape}]
\fill (1,1) circle (2pt) ;\fill (1,2) circle (2pt) ;\fill (1,3) circle (2pt) ;\fill (1,4) circle (2pt) ;\fill (2,1) circle (2pt) ;\fill (2,2) circle (2pt) ;\fill (2,3) circle (2pt) ;\fill (3,1) circle (2pt) ;\fill (3,2) circle (2pt) ;\fill (3,3) circle (2pt) ;\fill (4,1) circle (2pt) ;\fill (4,2) circle (2pt) ;\fill (4,3) circle (2pt) ;\fill (5,1) circle (2pt) ;\fill (5,2) circle (2pt) ;\fill (5,3) circle (2pt) ;\fill (6,0) circle (2pt) ;\fill (6,1) circle (2pt) ;\fill (6,2) circle (2pt) ;\fill (6,3) circle (2pt) ;\fill (8,0) circle (2pt) ;\fill (8,1) circle (2pt) ;\fill (8,2) circle (2pt) ;\fill (8,3) circle (2pt) ;\fill (8,4) circle (2pt) ;\fill (9,0) circle (2pt) ;\fill (9,1) circle (2pt) ;\fill (9,2) circle (2pt) ;\fill (9,3) circle (2pt) ;\fill (9,4) circle (2pt) ;\fill (10,0) circle (2pt) ;\fill (10,1) circle (2pt) ;\fill (10,2) circle (2pt) ;\fill (10,3) circle (2pt) ;\fill (10,4) circle (2pt) ;\fill (11,0) circle (2pt) ;\fill (11,1) circle (2pt) ;\fill (11,2) circle (2pt) ;\fill (11,3) circle (2pt) ;\fill (11,4) circle (2pt) ;\fill (12,0) circle (2pt) ;\fill (12,1) circle (2pt) ;\fill (12,2) circle (2pt) ;\fill (12,3) circle (2pt) ;\fill (12,4) circle (2pt) ;\fill (13,0) circle (2pt) ;\fill (13,1) circle (2pt) ;\fill (13,2) circle (2pt) ;\fill (13,3) circle (2pt) ;
\draw (6,0) -- (6,3) -- (5,3) -- (5,1) -- (4,1) -- (4,3) -- (3,3) -- (3,1) -- (1,1)      --  (1,2) -- (2,2) -- (2,3) -- (1,3) -- (1,4) ;
\draw (8,0) -- (8,4) -- (9,4) -- (9,0) -- (10,0) -- (10,4) -- (11,4) -- (11,0)        --  (13,0) -- (13,1) -- (12,1) -- (12,2) -- (13,2) -- (13,3) -- (12,3) -- (12,4)    ;
\node at (6.3, 0) {\tiny 0} ;  \node at (6.3, 1) {\tiny 6} ;\node at (6.3, 3) {\tiny 18} ; \node at (0.7, 1) {\tiny 1} ; \node at (0.7, 4) {\tiny 19} ; \node at (7.7, 0) {\tiny 0} ; \node at (7.7, 1) {\tiny 6} ;\node at (7.7, 4) {\tiny 24} ; \node at (13.3, 0) {\tiny 5} ; \node at (13.3, 3) {\tiny 23} ; \node at (12.3, 4) {\tiny 28} ; 
\end{tikzpicture}\end{center}
\end{figure}
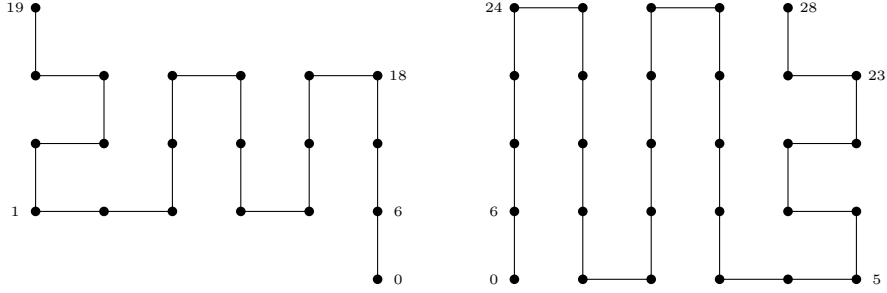

\begin{proof}
When~$x$ is even, the~$\omega$-construction with support~$\{1,x\}$ finishes on the smallest element of a fauxset.  For a perfect realization, we need to finish on the largest one.  We adapt the $\omega$-construction to achieve this by moving upwards through both of the last two fauxsets alternately, hopping back and forth with 1-edges.

First, let~$v = (x+2s-2) + (2sx-2s) + 1 = (2s+1)x -1$.  
Let~$\Theta_s$ be the path 
$$[x-2, x-1, 2x-1, 2x-2, 3x-2, 3x-1, \ldots ,
  (2s-1)x - 1 , 2sx-1 , 2sx - 2  , (2s+1)x -2 ].$$
This uses all of the elements in the fauxsets~$\phi_{x-2}$ and~$\phi_{x-1}$ and realizes~$\{1^{2s}, x^{2s}\}$.

The required linear realization is
$$\left( \apph{}{k=0}{x-3} \Psi_k \right) \uplus \Theta_s,$$
starting at~0.
The portion in parenthesis realizes~$\{1^{x-3} , x^{2s(x-2)} \}$ and so the path as a whole realizes
$$ \{1^{x-3} , x^{2s(x-2)} \} \cup \{ 1 \} \cup \{1^{2s}, x^{2s}\} = \{ 1^{x+2s-2} , x^{ 2sx -2s} \}.$$ 
The first element is~0 and the final element is~$(2s+1)x -2 = v-1$, so the realization is perfect.

The second construction adapts the $\omega$-construction~$\mathbf{h_2}$ in the analogous way.
Let $v = (x+2s-1) + (2sx - 2s + 2) + 1 = (2s+1)x + 2$.  Let~$\Phi_s$ be the path
$$[ 2, 1, x+1, x+2, 2x+2, 2x + 1, \ldots ,
     (2s-1)x+2  , 2sx + 2 , 2sx + 1 , (2s+1)x + 1 ].$$
The required linear realization is
$$ \Psi_0 \uplus \left( \apph{}{k=y-1}{3} \Psi_k \right) \uplus \Phi_s .$$
Checking that this is perfect and realizes the correct multiset is similar to the first case.

Figure~\ref{fig:even_perf} gives examples of the construction for~$\{1^7,6^{12}\}$ and~$\{1^8,6^{20}\}$.
\end{proof}

\begin{lem}\label{lem:midperf}{\rm (Reordering one generic fauxset)}
Suppose a standard realization for a multiset~$L$ includes a subsequence of the form~$\Psi_k^{(s_1,s_2)}$, a fauxset segment with respect to~$x$ in order with~$w = s_2 - s_1$ edges.  Let~$M$ be a multiset of size~$w$ that has a perfect realization.  Then~$(L \setminus \{x^w\}) \cup xM$ has a standard realization.
\end{lem}

\begin{proof}
Let $\Psi_k^{(s_1,s_2)} = [s_1x+k, (s_1+1)x + k, \ldots s_2x+k]$, which realizes~$\{x^w\}$, and the perfect realization for~$M$ be $\mathbf{h} = [h_1, h_2, \ldots, h_{s+1}]$.  Then
$$(\mathbf{h}+s_1)x + k = [(h_1+s)x + k, (h_2+s)x + k, \ldots, (h_{s+1}x+s) + k]$$
realizes~$xM$, starts and finishes on the same vertices as~$\Psi_k^{(s_1,s_2)}$ and uses the same vertices as~$\Psi_k^{(s_1,s_2)}$.

Replace~$\Psi_k^{(s_1,s_2)}$ with~$(\mathbf{h}+s_1)x+k$ to obtain the desired realization.
\end{proof}

We can now prove the second main result of the section.

\begin{thm}\label{th:1xtx_tight}{\rm (Odd $x$ improved)}
Let~$L = \{1^a, x^b, tx^c \}$ be a multiset of size~$v-1$ with~$x$ odd and~$v \geq (t+2)x$.  Write $v = qx + r$ with $0 \leq r < x$ and let~$s = \lfloor q/2t \rfloor$.    If~$a \geq \omega(x,b+c)$ then~$L$ has a standard linear realization when any one of the following conditions hold:
\begin{enumerate}
\item  $r$ is even and $b \geq tx$,
\item $r$ and $t$ are odd and $b \geq tx + 2t - 3$,
\item $r$ is odd, $t$ is even and $b \geq 4s + tx + 4t -9$,
\end{enumerate}
\end{thm}

\begin{proof}
In all cases we use use the basic shape of the argument of the proof of Theorem~\ref{th:1xtx_gen} and indicate how augmenting it in particular ways leads to the desired bounds.

Consider Case 1.  We use the $\omega$-construction~$\mathbf{h_1}$.  We may pair all fauxsets except the final one and we do not have the dangling edge requiring the~``+1".  So the bound on~$b$ becomes
$$ b \geq (x-1)t + t  = tx.$$

Consider Case 2.   We use the $\omega$-construction~$\mathbf{h_3}$.  We may use Lemma~\ref{lem:midperf} on each of the two unaltered fauxsets.  The $\omega$-construction with support~$\{1,t\}$ is perfect when the number of vertices is a multiple of~$t$.  If the fauxsets have~$q^*$ vertices then we may use a perfect realization with $\lfloor q^*/t \rfloor t$ vertices (requiring~$t-1$ 1-edges) and append up to~$t-1$ further 1-edges to reach the required length.  The bound on~$b$ becomes
$$ b \geq (x-3)t + 1 + t + 2\cdot2(t-1) = tx+2t-3.$$

Consider Case 3.  The idea is the same as Case~2, but as~$t$ is even we use a perfect realization with support~$\{1,t\}$ from Lemma~\ref{lem:even_perf}.  

The fauxsets in question have at least~$q$ vertices. 
There is a perfect realization  with support~$\{1,t\}$ whenever the number of vertices is congruent to $(t-1), (t+2) \pmod{2t}$.  This means that at most $(2t-3)$ 1-edges are required to expand this to a perfect linear realization with~$q$ vertices by concatenating a perfect linear realization with support~$\{1\}$ of an appropriate length.  The realization itself uses at most $(t+2s-2)$ 1-edges,   so $(2t-3) + (t+2s-2) = (3t+2s-5)$ 1-edges are sufficient for a perfect linear realization with~$q$ vertices to be guaranteed to exist.
The bound on~$b$ becomes
$$b \geq (x-3)t + 1 + t + 2(3t+2s-5) = 4s + tx + 4t -9. $$
\end{proof}

Observe that in the first two cases of Theorem~\ref{th:1xtx_tight} (which apply unless both~$r$ is odd and~$t$ is even) the proportion $b/(v-1)$ of $x$-edges required tends to~0 as~$v$ increases.
In the remaining case, the proportion of~$x$-edges improves by a factor of~$t$ compared to Theorem~\ref{th:1xtx_gen} as $v$ increases, to $2/tx = 2/y$.

\bigskip

By limiting the congruence classes of~$v$ under consideration we are able to use a different approach to prove a general result for multisets of the form~$\{1,x,tx\}$.

\begin{thm}\label{th:1xtx}{\rm (General case)}
Let~$x \geq 3$ and~$t \geq 2$ with~$(t,x) \neq (2,3)$. Suppose~$v \equiv \pm 1 \pmod{tx}$.  Then multisets of size~$v-1$ with support~$\{ 1, x, tx\}$ are realizable provided
$$  v > \frac{t^2x^2 + t^2x + tx^2  + t + x + 1}{tx-t-x-1}.$$
In particular, such multisets are realizable when~$t,x \geq 7$.
\end{thm}

\begin{proof}
Let $L = \{ 1^a, x^b, tx^c\}$ and suppose~$v = ktx \pm 1$.  Necessarily, $k \geq 2$ as if~$k=1$ we have~$\widehat{tx} = 1$.  

As~$x^{-1} = \pm tk \pmod{v}$, we have 
$$L' = \{ \widehat{tk}^a, 1^b, \widehat{kt^2x}^c \} = \{ 1^b, kt^a, t^c \}.$$
As~$(tx)^{-1} = \pm k \pmod{v}$, we have 
$$L'' = \{ \widehat{k}^a, \widehat{kx}^b, 1^c \} = \{ 1^c, k^a, kx^b \}.$$
Therefore, by Theorem~\ref{th:1xy_new}, any counterexample has $a < x+tx-1$, $b < kt + t-1$ and $c < k + kx-1$ and so satisfies
$$a+b+c < x+tx+kt+t+k+kx -1 = k(t+x+1) + x+t+tx-1.$$
As $a+b+c=v-1 = (ktx \pm 1) -1$, a counterexample satisfies
$$ktx < k(t+x+1) +x+t+tx+1.$$
Equivalently,
$$k < (tx + t +x  +1) / (tx - t - x -1).$$

Given fixed~$t \geq 2$ and~$x \geq 3$ (and noting that the denominator only vanishes in the excluded case~$(t,x) = (2,3)$) this inequality can only hold for small~$k$.  Combining this with the inequality~$v \geq ktx-1$ gives the main result.

Suppose~$t,x \geq 7$.  We consider when the above inequality implies that~$k < 2$, a contradiction.  Re-arranging, we find that
$$(tx + t +x  +1) / (tx - t - x -1) < 2$$
if and only if $3t + 3x + 3 < tx$.  This is true for~$t=x=7$, as~$45 < 49$, and increasing either~$x$ or~$t$ (or both) maintains the truth of the inequality.  Hence the subsidiary statement holds.
\end{proof}

As $v \equiv \pm 1 \pmod{tx}$ in Theorem~\ref{th:1xtx}, all the cases covered fall under the Coprime BHR~Conjecture.  Also note that the excluded case~$(t,x) = (2,3)$ corresponds to the support~$\{1,3,6\}$, for which the BHR~Conjecture is known to hold~\cite{CO}.

\section{Some consequences for the BHR Conjectures}\label{sec:conseq}

What else do the approaches and tools developed here have to say about the three conjectures?  We first consider the support~$\{1,3,y\}$ and then see how the different tools come together to make progress with the Coprime BHR Conjecture for the supports~$\{1, 6, 18\}$.

In a recent paper, Avila completes the first work focused on supports of the form $\{1,3,y\}$, proving that admissible multisets~$\{1^a,3^b,y^c\}$ are realizable when~$y$ is even, $c$ is odd and either~$a \geq y+1$ or~$a=y$ and~$3 \nmid b$~\cite{Avila23}.  From Theorem~\ref{th:1xy2} we know that such multisets are realizable when~$a \geq y+3$ and when~$y$ is even and~$a = y+2$.  From Theorem~\ref{th:known} we know they are realizable when~$v \leq 37$ and when~$y \leq 7$.

Theorems~\ref{th:13y} and~\ref{th:13ynum} demonstrate how the tools developed in this paper improve on these results. 

\begin{thm}\label{th:13y}{\rm (Applications to $x=3$)}
Let~$L = \{1^a, 3^b, y^c\}$ be an admissible multiset of size~$v-1$. Write~$v = 3q + r$ with $0 \leq r < 3$.  Then~$L$ has a linear realization in the following cases:
\begin{enumerate}
\item $a \geq y+1$,
\item $a \geq 7$ and~$b \geq y-4$, 
\item $a \geq 4$,  $a+b \geq y+1$ and $c \equiv 0,1 \pmod{y}$,
\item $a \geq 2$, $y=3t$, $v \geq 3t+6$  and either:
\begin{itemize}
\item $r$ is even and~$b \geq 3t$,
\item $r$ is odd and $b \geq 3t+1$.
\end{itemize}
\end{enumerate}
\end{thm}

\begin{proof}
Item~1 follows from Theorem~\ref{th:1xy_new}.

Item~2 follows, with a little more work, from Theorem~\ref{th:odd_hops} and its proof.    
Applying that result as stated to the case~$x=3$ shows that it is sufficient to take $a \geq 9$ and~$b \geq y-4$.
Tracking the proof of Theorem~\ref{th:odd_hops} with the extra knowledge that $x=3$, we find two spots where we may make improvements.  

In the first realization, we assume that~$x$ edges are required, as for general odd~$x$ it might be that~$\omega(x,d) = x$ for some~$d$ under consideration.  When~$x=3$ we have $\omega(x,d) = x-1 = 2$ for all~$d$ and we may reduce the number of 1-edges required by~1.  The same reasoning applies when considering the standard linear realization with~$v_2$ vertices, allowing us to reduce the number of 1-edges by 1 again.
We conclude that it is sufficient to take~$a \geq 4x-3 -2 = 7$ when~$x=3$.

Item~3 follows from Theorem~\ref{th:thicken}.    
 
Item~4 follows from Theorem~\ref{th:1xtx_tight}.1; the first bullet immediately and the second after an additional observation.
If $r$ is odd then necessarily~$r=1$.  In this case we may follow Case~1 of Theorem~\ref{th:1xtx_tight} with the $\omega$-construction~$\mathbf{h_2}$.  The only difference to the argument is that we do have the ``dangling edge requiring the +1" and so we obtain a bound of $tx+1 = 3t+1$ rather than~$tx$.   \end{proof}

\begin{thm}\label{th:13ynum}{\rm (Another application to $x=3$)}
The  BHR Conjecture holds for $\{1,3,3t\}$ when $v \equiv \pm 1 \pmod{3t}$ and $v > 6t+53$.
\end{thm}

\begin{proof}
The BHR Conjecture holds for~$\{1,3,6\}$~\cite{CO}, so let~$t \geq 3$.
Applying Theorem~\ref{th:1xtx}, we find that the multiset is guaranteed to be realizable when
$$v >  \frac{12t^2 + 19t + 4}{2t-4}.$$
The result follows upon noting that
$$\frac{12t^2 + 19t + 4}{2t-4} =  6t + 17 +   \frac{36}{t-2} \leq 6t+53.$$
\end{proof}

Theorem~\ref{th:1618} illustrates what may happen when using the old and new tools to consider a fixed support, in this case~$\{1, 6, 18\}$.  Rather than aiming for efficiency, the proof follows the natural steps when faced with a question like this: use the most general (and, usually, easiest to apply) tool first and then successively apply more specialized complex tools to chip away at the remaining possible counterexamples.

\begin{thm}\label{th:1618}{\rm (Fixed support application example)}
The Coprime BHR Conjecture holds for the support~$\{1,6,18\}$, except possibly for~$3$ values of~$v$: $47, 49, 59$.
\end{thm}

\begin{proof}
Let~$v = 6k \pm 1$.   Then~$\widehat{6^{-1}} = k$, so $L' = \{1,3,k\}$.  Also,~$\widehat{3^{-1}} = 2k$, so $L'' = \{1, 2k, \widehat{3^{-1}k}\}$.
By Corollary~\ref{cor:1xy}, any counterexample must have
$$ 6 + 18 + 3 + k + 2k + \widehat{3^{-1}k}  = 3k+27 + \widehat{3^{-1}k} \geq v+2$$
As~$v \not\equiv 0 \pmod{3}$, it must be that~3 divides one of~$k$, $v+k$ and~$2v+k$.  If $3\mid k$ then $\widehat{3^{-1}k} = k/3$; if $3 \mid v+k$ then $\widehat{3^{-1}k} = (7k \pm 1)/3$; if  $3 \mid 2v+k$ then $\widehat{3^{-1}k} = \widehat{(13k \pm 2)/3} = (5k \pm 2)/3$.  So  $\widehat{3^{-1}k} \leq (7k + 1)/3$.

Therefore, we cannot have a counterexample when
$3k+27 +(7k+1)/3 < 6k+1$.  Simplifying, there is no counterexample when~$k \geq 40$; that is $v \geq 6\cdot 40 -1 = 239$.

The remainder of the proof uses a computer to assess the arithmetic conditions and keep track of cases, but not to construct realizations.  It turns out that the method above also covers all but~32 of the cases of the Coprime BHR Conjecture for~$\{1,6,18\}$ with~$v \leq 239$.  
We shorten this list with previously known and new results.  

As~$\{1,3,k\}$ is realizable for~$k \leq 7$ by Theorem~\ref{th:known}.2, any counterexample must have $v \geq 6\cdot8 -1 = 47$, eliminating three of the remaining cases and verifying the lower bound on~$v$ in the statement of the theorem.  

Using Theorem~\ref{th:1xy_new} rather than Corollary~\ref{cor:1xy}  eliminates a further~21 cases, reducing the total to~8.  The largest remaining value of~$v$ is~131.

Consider~$v=131$.  The three equivalent multisets are $\{1^a,6^b,18^c\}$, $\{1^b, 3^c, 22^a\}$ and~$\{1^c, 44^b, 51^c\}$.  By Theorem~\ref{th:1xy_new}, any counterexample must have~$a \leq 16$, $b \leq 23$ and~$c \leq 92$.  Hence we also have $b = (v-1) - (a+c) \geq 16$ and $c = (v-1) - (a+b) \geq 91$. We may therefore apply Theorem~\ref{th:13y}.2 to $L' = \{1^b, 3^c, 22^a\}$ to rule out all possible counterexamples. The same approach works for~$v \in \{ 61, 85, 95 , 121\}$, 
leaving~3 possible values of~$v$ at which a counterexample has not been ruled out: $v \in \{47, 49, 59 \}$. 
\end{proof}



Although the methods used in the proof of Theorem~\ref{th:1618} do not cover those last~3 values of~$v$, in each case they do significantly limit the form a counterexample might take.    
For example, only~15 possible counterexamples remain for $v = 59$.  A further~5 of these may be eliminated by applying Theorem~\ref{th:xoddyeven}, leaving the following~10 possibilities for~$(a,b,c)$:
$$ ( 13, 6, 39 ), ( 14, 5, 39 ), ( 14, 6, 38 ), ( 15, 4, 39 ), ( 15, 5, 38 ), $$ $$ ( 15, 6, 37 ), ( 16, 3, 39 ), ( 16, 4, 38 ), ( 16, 5, 37 ), 
  ( 16, 6, 36 ). $$


As a final illustration of the approach we show that for supports~$\{1,x,y\}$ with~$y$ sufficiently large compared to~$x$ there are infinitely many values of~$v$ at which the BHR~Conjecture holds. 

\begin{thm}\label{th:bigy}{\rm (A general result)}
Let~$\supp(L) = \{1,x,y\}$ with~$x>2$ and $y > (2x^2 + 2x + 1) / (x-2)$.   For all sufficiently large~$v$ with~$v \equiv \pm 1 \pmod{xy}$, the BHR Conjecture holds for~$L$.
\end{thm}

\begin{proof}
Let~$v = xyk \pm 1$ and let~$L = \{1^a, x^b, y^c\}$.  We have~$\widehat{x^{-1}} = yk$ and so~$L' = \{1^b, (yk)^a, \widehat{y^2k}^c\}$, and  $\widehat{y^{-1}} = xk$ and so~$L'' = \{1^c, (xk)^a, \widehat{x^2k}^c\} = \{1^c, (xk)^a, (x^2k)^c\}$.  In the notation of Theorem~\ref{th:1xy_new} we have
\begin{eqnarray*}
 f(x,y) + + f(\dot{x} ,\dot{y}) + f(\ddot{x} , \ddot{y})
   &\leq& (x+y-1) + (xk+x^2k-1) + (yk+\widehat{y^2k} -1) \\
   & = & (x+y-3) + (x+x^2+y)k + \widehat{y^2k} \\
   & \leq & (x+y-3 + \frac{1}{2}) + (x+x^2+y + \frac{xy}{2})k, 
\end{eqnarray*}
as $\widehat{y^2k} \leq v/2 \leq (xyk+1)/2$.

If there is an~$\epsilon > 0$ with $(x+x^2+y + \frac{xy}{2}) < xy - \epsilon$ then we obtain
$$ f(x,y) + + f(\dot{x} ,\dot{y}) + f(\ddot{x} , \ddot{y}) 
 < (x+y-\frac{5}{2}) + (xy - \epsilon)k$$ 
which is less than~$v+2$ for sufficiently large~$k$, and therefore for sufficiently large~$v$.
Set~$\epsilon = 1/2$.  Then $(x+x^2+y + \frac{xy}{2}) < xy - \epsilon$ exactly when $y > (2x^2 + 2x + 1) / (x-2)$ and the result follows.
\end{proof}

As 
$$\frac{2x^2 + 2x + 1}{x-2} = 2x + 6 + \frac{13}{x-2} \leq 2x+19$$
for~$x>2$, at the expense of a slight weakening of the result we can simplify the bound on~$y$ in Theorem~\ref{th:bigy} to~$y > 2x+19$.

We have~$\gcd(xyk \pm 1, x) = 1 = \gcd(xyk \pm 1, y)$, hence all of the results generated by Theorem~\ref{th:bigy} are instances of the Coprime BHR Conjecture.  As~$\gcd(xy, \pm1) =1$, by the Dirichlet prime number theorem there are infinitely many primes of the form $xyk \pm 1$.  Therefore, the theorem implies that Buratti's Conjecture holds for infinitely many supports with the same conditions.

\end{document}